\documentclass[oneside,11pt]{amsart}
\usepackage{mathrsfs}
\usepackage{mathdg2015}
\usepackage{geometry}
\usepackage{graphicx}
\usepackage[new]{old-arrows}
 \DeclareRobustCommand\longtwoheadrightarrow
     {\relbar\joinrel\twoheadrightarrow}
\usepackage{tikz-cd}
\usepackage{textgreek}
\usepackage[hypertexnames=false,pdfencoding=auto]{hyperref}

\expandafter\def\csname ver@etex.sty\endcsname{3000/12/31}

\usepackage{autonum}

\newcommand{\p}{\prime}
\renewcommand{\C}{\mathbb{C}}

\newcommand{\Dnu}{D_\nu}
\newcommand{\textBC}{{\upshape\textsc{bc}}}
\newcommand{\BC}{\ensuremath{\text{\resizebox{!}{4pt}{\tiny\upshape BC}}}}
\newcommand{\s}{\mathrm{s}}
\newcommand{\dsX}{{\partial_{\s}X}}
\newcommand{\dsiX}{{\partial_{\s,i}X}}
\newcommand{\dsM}{{\partial_{\s}M}}
\newcommand{\dBCX}{{\partial_\BC X}}
\newcommand{\dBCM}{{\partial_\BC M}}
\newcommand{\dsBCX}{{\partial_{\s,\BC} X}}
\newcommand{\dsBCiX}{{\partial_{\s,\BC,i} X}}

\renewcommand{\dX}{{\partial X}}
\newcommand{\dM}{{\partial M}}

\DeclareMathOperator{\spn}{span}

\newcommand{\codim}{\operatorname{codim}}

\newcommand{\phisigma}{{}^\phi\sigma} 

\newcommand{\phiT}{{}^\phi T} 
\renewcommand{\phihat}{{\hat{\phi}}}
\newcommand{\phiBC}{{\phi_{\BC}}}
\newcommand{\phiBCi}{{\phi_{\BC,i}}}
\newcommand{\sigmafull}{{}^\phi\sigma_{\mathrm{full}}}
\newcommand{\DN}{\mathcal{N}}

\newcommand{\bbf}{\mathrm{bf}}

\newcommand{\SymP}{\mathfrak{p}}
\newcommand{\Jp}{\hspace{-2pt}\mathscr{J}}
 \newcommand{\kerningJp}{\hspace{-2pt}\mathscr{J}}
\newcommand{\Jphat}{\hspace{3pt}\hat{\hspace{-3pt}\Jp}}
 \newcommand{\kerningJphat}{\hspace{1pt}\hat{\hspace{-3pt}\Jp}}
\newcommand{\Cdotinf}{{\dot C^\infty}}
 \newcommand{\CdotBCinf}{{\dot C^\infty_\BC}}
 \newcommand{\Cdotsinf}{{\dot C^\infty_{\s}}}
 \newcommand{\CdotsBCinf}{{\dot C^\infty_{\s,\BC}}}
\newcommand{\calFdotBC}{\dot{\calF}_\BC}
\newcommand{\Phatbar}{{\bar\Phat}}
\newcommand{\Sph}{\mathbb{S}}
\newcommand{\B}{\mathbb{B}}
\newcommand{\cl}{{\mathrm{cl}}}
\newcommand{\Vcomp}{V_{\mathrm{comp}}}
\newcommand{\Picomp}{\Pi_{\mathrm{comp}}}
\newcommand{\Vsh}{V_{\mathrm{sh}}}
\newcommand{\Pish}{\Pi_{\mathrm{sh}}}
\newtheorem{assumption}[theorem]{Assumption}

\date{\today}
\title[Calder\'on Projector for Fibred Cusps]{The Calder\'on Projector for Fibred Cusp Operators}
\author{Karsten Fritzsch}
\address{Institut f\"ur Analysis, Leibniz Universit\"at Hannover}
\email{k.fritzsch@math.uni-hannover.de}
\author{Daniel Grieser}
\address{Institut f\"ur Mathematik, Carl von Ossietzky Universit\"at Oldenburg}
\email{daniel.grieser@uol.de}
\thanks{Part of this work was done while DG was in residence at the Mathematical Sciences Research Institute in Berkeley, California, during the Fall semester 2019, supported by the National Science Foundation under Grant No. DMV-1440140.
DG and ES were partially supported by DFG Priority Programme 2026 \lq Geometry at Infinity\rq.}
\author{Elmar Schrohe}
\address{Institut f\"ur Analysis, Leibniz Universit\"at Hannover}
\email{schrohe@math.uni-hannover.de}
\begin{document}

\tikzcdset{row sep/normal=1em}



\begin{abstract}
A Calder\'on projector for an elliptic operator $P$ on a manifold with boundary $X$ is a projection from general boundary data to the set of boundary data of solutions $u$ of $Pu=0$. Seeley proved in 1966 that for compact $X$ and for $P$ uniformly elliptic up to the boundary there is a Calder\'on projector which is a pseudodifferential operator on $\partial X$. We generalize this result to the setting of fibred cusp operators, a class of elliptic operators on certain non-compact manifolds having a special fibred structure at infinity.
This applies, for example, to the Laplacian on certain locally symmetric spaces
or on particular singular spaces, such as a domain with cusp singularity or the complement of two touching smooth strictly convex domains in Euclidean space. Our main technical tool is the $\phi$-pseudodifferential calculus introduced by Mazzeo and Melrose.

In our presentation we provide a setting that may  be useful for doing analogous constructions for other types of singularities.
\end{abstract}

\maketitle

\tableofcontents

\section{Introduction}

A Calder\'on projector associated with an  elliptic partial differential operator $P$ of order $m$ on a compact manifold  $X$ with non-empty boundary $\dX$ is a projection $C$ in $\Cinf(\partial X)^m$ to the set of boundary data of solutions of the homogeneous equation
\begin{equation}
 \label{eqn:cauchy data space}
 \{(u_{|\partial X},\Dnu u_{|\partial X},\dots, \Dnu^{m-1} u_{|\partial X}):\, u\in\Cinf(X), Pu = 0\} \,.
\end{equation}
Here $\nu$ is some choice of transversal vector field in a neighborhood of the boundary and $\Dnu = \frac{1}{i}\partial_\nu$.

It was first observed by Calder\'on \cite{Cal63} that such a projection exists which is a pseudodifferential operator, with an explicit principal symbol, and that this can be used to study boundary value problems. The first complete proof was given by Seeley in 1966, \cite{See66,See69}.
The result applies more generally to operators acting between sections of vector bundles.

In the present paper we extend this result to certain settings where the boundary has singularities, that is, we construct Calder\'on projectors with full control of their behavior near the singularities. We consider certain classes of singularities, often called of fibred cusp type, but our purpose is also to present this construction in a systematic way to lay the foundation for generalizing the analysis to other types of singularities.
Typical examples to which our results apply are the Laplace-Beltrami or Dirac type operators on  Riemannian manifolds with boundary which locally are of one of the following types:
\begin{enumerate}
 \item
 Domains with (incomplete) cusp singularity such as $\{(\xi,\eta)\in\R^m\times\R^k: |\xi|\leq |\eta|^2\}$ near $(\xi,\eta)=0$ where $m,k\in\N$; for $k=1$ this is what is commonly called an incomplete cusp; for $m=1$ this has the same geometry as the complement of two touching strictly convex sets in $\R^{k+1}$ as in Figure \ref{fig:blow-up example} left;
 \item
 certain types of
 domains in locally symmetric spaces of $\Q$-rank one, for example the strip given by $\{|\Re w| \leq \frac14,\ \Im w \geq 1\}$ in the complex plane with the hyperbolic metric (considered locally near infinity),
 \item
 spaces of the form $C\times F$, where $C$ is the \lq far end\rq\ of a cone (e.g.\ the outside of a ball in $\R^n$) and $F$ is a compact Riemannian manifold with non-empty boundary,
 or more generally fibre-bundles of this sort.
\end{enumerate}
In the second and third example the \lq singularity\rq\ of the boundary is at infinity, so here we give a uniform description of the behavior of the Calder\'on projector at infinity. We also allow spaces to have several such singularities of different types, as well as analogous types of singularities away from the boundary, for example the exterior of a smooth bounded domain in $\R^n$. More details on these examples are provided in Section \ref{ssec:ex}, and the precise class of spaces and operators is described below.

The main motivation for the study of Calder\'on projectors is their use in the analysis of regularity and Fredholm properties of boundary value problems.
In particular, if $P$ is the Laplace-Beltrami operator $\Delta$ associated to some Riemannian metric on a compact manifold $X$ with boundary,  then Seeley's result implies, by standard pseudodifferential calculus techniques, that the Dirichlet-Neumann operator, defined by
\begin{equation}
\label{eqn:def DN op}
 \DN: \Cinf(\partial X)\to\Cinf(\partial X),\ f \mapsto \partial_\nu u\ \text{ where } \ \Delta u = 0, u_{|\partial X} = f \,,
\end{equation}
with the outward unit normal derivative $\partial_\nu$, is a pseudodifferential operator of order 1, with principal symbol $|\xi|$.

We will apply our results on Calder\'on projectors to boundary value problems and the Dirichlet-Neumann operator in the fibred cusp setting in a separate paper.
The Dirichlet-Neumann operator occurs in a number of contexts, such as the Calder\'on inverse problem or the theory of water waves. It also appears in the context of the plasmonic  eigenvalue problem, \cite{Gri14}, and here the fibred cusp geometry is of particular importance.
This two-sided boundary value problem describes the coupling of electromagnetic fields to the electron gas of a conducting body. Since the geometry of the body can be used to specifically tailor properties of the resulting surface waves, its solutions on more singular spaces have seen great interest in recent years. See for instance \cite{GR09,GUB+09,Sav12,CCN18,Sch18} and the references therein.
If this body consists of two (nano-meter sized) balls that touch each other then we are led precisely to the geometry studied in this paper. This is also considered in \cite{Sch18} and the first author's thesis \cite{Fri14}. In the latter, the method of layer potentials is studied in the light of manifolds with corners and conormal distributions. See also \cite{Fri19} for the geometrically simpler case of the half-space.

\subsection{Main theorems} \label{subsec:main theorems}
The differential operators we consider have, in suitable local coordinates near the singularities, the form
\begin{equation}
 \label{eqn:phi-op local coords}
 P = x^{-cm} \sum_{k+|\alpha|+|\beta|\leq m} a_{k\alpha\beta}(x,y,z) (x^2 D_x)^k (x D_y)^\alpha D_z^\beta \,.
\end{equation}
Here $x>0$, $y=(y_1,\dots,y_b)\in U\subset\R^b$, $z=(z_1,\dots,z_f)\in V$ where $V\subset \R^f$ or $V\subset[0,\infty)\times\R^{f-1}$ with $b\in\N_0$, $f\in\N$ and open subsets $U$, $V$. As usual, we let $D_x=\frac1i\partial_x$ etc. The singularity corresponds to $x=0$, and the case  $V\subset[0,\infty)\times\R^{f-1}$ corresponds to a neighborhood of the boundary, which then is $z_1=0$. Also, $m\in\N$ is the order of $P$ and
$c\in\Z$ is a parameter. The coefficients $a_{k\alpha\beta}$ are assumed to be smooth up to $x=0$, and $P$ is assumed to be elliptic, in a uniform sense as $x\to0$ to be described later.
The singular structure of the operator is reflected in the $x$-factors in $x^2D_x$ and $xD_y$, and to a lesser extent in the prefactor $x^{-cm}$.

\subsubsection*{{\sc BC}-boundary and singular boundary} Globally it is useful to describe $P$ as acting on a compact manifold with corners $X$ which has two types of boundary hypersurfaces: the `singular boundary' $\dsX$, given by $x=0$, and the `boundary at which boundary conditions could be imposed', the `\textBC-boundary' denoted by $\dBCX$ and given locally by $z_1=0$. These intersect in their common boundary $\dsBCX$. Note that $P$ acts on functions defined on $X\setminus\dsX$ only, but adding the boundary $\dsX$ to the space allows us to express the singular behavior of $P$ and other objects efficiently.
The difference in the roles of the $y$ and $z$ variables in  \eqref{eqn:phi-op local coords} is geometrically and globally described by an additional piece of data: a fibration of the singular boundary
$$F - \dsX \overset{\phi}\to B$$
with compact base $B$ and compact fibre $F$. Locally, $\phi$ is just the map $\phi:(y,z)\mapsto y$, so $y$ are base coordinates and $z$ are fibre coordinates. The fibre $F$ has boundary given locally by $z_1=0$.
The restriction of $\phi$ to $\partial(\dsX) = \dsBCX$ is again denoted by $\phi$  and defines a fibration $\partial F - \dsBCX \overset{\,\phi\,}\to B$.
We call such a space $X$ a $\phi$-\textBC-\textbf{manifold} or $\phi$-\textbf{manifold with} \textBC-\textbf{boundary} and an operator as in \eqref{eqn:phi-op local coords}, with $c=0$, a \textbf{$\phi$-differential operator} or short \textbf{$\phi$-operator}.
The definition extends in a straightforward way to operators acting between sections of  vector bundles $E,E'$ over $X$. The class of these operators of order $m$ is denoted
$$\Diff_\phi^m(X;E,E') \,.$$
See Section \ref{sec:setting} and Appendix \ref{sec:app mwc} for basics on manifolds with corners and more details, including how the examples above fit into this framework.

\subsubsection*{The $\phi$-calculus} In the case that $X$ has empty \textBC-boundary, i.e.\ if the fibres of $\phi$ do not have a boundary, the class of $\phi$-operators was introduced and studied by Mazzeo and Melrose \cite{MM98} (see also \cite{GH09,GH14}).
They showed that interior regularity extends to conormal boundary regularity
if the operator $P$ is \textbf{$\phi$-elliptic}, i.e.\ its $\phi$-principal symbol $\phisigma_m(P)$ is invertible, and that $P$ is Fredholm in naturally associated $L^2$ spaces if and only if it satisfies the stronger condition of being \textbf{fully elliptic}, which in addition requires invertibility of a family $N(C) (\mu)$, $ \mu \in \R \times T^*B$, of differential operators, called \textbf{normal family} of $P$, on the fibres of $\phi$.
These results are proved via construction of a pseudodifferential operator ($\Psi$DO) calculus adapted to the singular structure.
The space of such $\phi$-$\Psi$DOs of order $m\in\R$ is denoted $\Psi^m_\phi(X)$.
See Appendix \ref{sec:app phi ops} for more details. 

In the setting of $\phi$-\textBC-manifolds with non-empty \textBC-boundary, $\dBCX$ with the restricted fibration is a $\phi$-manifold (without \textBC-boundary), so $\Psi^*_\phi(\dBCX)$, the space of $\phi$-$\Psi$DOs with respect to the fibration $\phi : \partial(\dBCX) = \dsBCX \rightarrow B$, is defined.

\subsubsection*{Function spaces} When considering Calder\'on projectors on a singular or non-compact space we need to specify the growth behavior at the singular set which we allow for the functions $u$ and for the vector field $\nu$ in \eqref{eqn:cauchy data space}.
The geometrically natural condition on $\nu$ is that it is a $\phi$-vector field, see \eqref{eqn:def phi-vector field}, transversal to $\dBCX$.
For the functions $u$ there are various ways to restrict the growth behavior at $\dsX$. We use the letter $\calF$ to denote any choice of function space encoding such behavior.
Thus, we have spaces $\calF(X)$, $\calF(\dBCX)$ of functions on $X$ and $\dBCX$, and similarly spaces of sections of vector bundles.
We call $\calF$ \textbf{admissible} if it behaves well under restriction to $\dBCX$ and under the action of $\phi$-$\Psi$DOs, see Definition \ref{def:admissible function space}. Also, we require all functions to be smooth in $X\setminus\dsX$ (see below for $L^2$-based spaces).
A simple example is $\calF(X) = x^\alpha \Cinf(X)$, where $\alpha\in\R$. This would allow $x^{\alpha}$ decay (if $\alpha>0$) or growth (if $\alpha<0$) at the singular set.
Other choices include spaces characterized by $L^2$- or $L^\infty$-based bounds (conormality) or full asymptotic expansions (polyhomogeneity) at $\dsX$, introduced in Appendix \ref{sec:app mwc}. The smallest choice of $\calF$ is $\Cdotsinf$, the space of functions vanishing to infinite order at $\dsX$ (in the geometric setting (iii)
this corresponds to functions rapidly decreasing at infinity with all derivatives), and the largest choice is $\calA_\s$, the space of functions conormal at $\dsX$.

\subsubsection*{Boundary data and Calder\'on projectors} We fix a $\phi$-vector field $\nu$ on $X$ transversal to $\dBCX$ throughout the paper. Then for any $m \in \N$ we have the boundary data map for admissible $\calF$
 \begin{equation}\label{eqn:first.gamma}
  \gamma : \calF(X) \longrightarrow \calF(\dBCX)^m \,,\quad
  u \longmapsto \gamma u = \big(u_{|\dBCX},\Dnu u_{|\dBCX}, \dotsc, \Dnu^{m-1}u_{|\dBCX}\big) \,.
 \end{equation}
In order to define $\gamma$ on sections of a bundle $E$, we also need to choose a connection on $E$ if $m>1$. We denote the \textbf{$\calF$-boundary data space} of $P \in \Diff_\phi^m(X;E,E')$ by
 \begin{equation}\label{eqn:first.bdy.data}
  \calB_{P,\calF} :=\{\gamma u :\, u\in\calF(X;E),\ Pu=0\} \subset \calF(\dBCX;E)^m
 \end{equation}
and we define an  \textbf{$\calF$-Calder\'on projector} for $P$ to be a projection in $\calF(\dBCX;E)^m$ to $\calB_{P,\calF}$. Note that there are many such projections.

A particular issue in dealing with the Calder\'on projector is that $P$ may have \textbf{shadow solutions}, i.e.\ sections $u\not\equiv0$ satisfying $Pu=0$ whose boundary data at $\dBCX$ vanish. This is a type of failure of unique continuation for $P$.
While this is not a problem for $P$ itself, we need to exclude it for the normal families of $P$ and its adjoint $P^\star$. A precise formulation is given in Assumption \ref{UCNF}.
Our main result is:
\begin{theorem}
 \label{thm:Calderon}
 Let $X$ be a $\phi$-manifold with non-empty \textBC-boundary and $E,E'$ be complex vector bundles over $X$. Let $c\in\Z$ and
 $$P =x^{-cm}\tilde P,\ \tilde P\in\Diff_\phi^m(X;E,E') \,,$$
 where $\tilde P$ is $\phi$-elliptic and satisfies Assumption \ref{UCNF}. Then:
\begin{enumerate}
\item There is an operator $C\in\Psi^*_\phi(\dBCX;E^m)$ which for any choice of admissible function space $\calF$  is an $\calF$-Calder\'on projector for $P$.
\\
When considered as acting between $m$-tuples of sections of $E$, the operator $C$ is an $m\times m$ matrix $(C_{kl})_{k,l=1\dots m}$,
where $C_{kl} \in \Psi_\phi^{k-l}(\dBCX;E)$.
\item   The $\phi$-principal symbol  $\phisigma_*(C)(p,\xi^\p)\in \End(E_p^m)$ of $C$ at any $(p,\xi^\p) \in \phiT^*\dBCX\setminus \{0\}$ is a Calder\'on projector 
for the differential operator $\phisigma_m(P)(p;D_t,\xi^\p)$ induced by the $\phi$-principal symbol $\phisigma_m(P)$ of $P$. More precisely, its range is the space of boundary data at $t=0$ of the solutions $v:\R\to E_p$ of \\ \mbox{$\phisigma_m(P)(p;D_t,\xi^\p)v(t)=0$} that vanish to infinite order at $t=+\infty$. 
\item The normal operator $N(C) (\mu)$ of $C$ at any  $ \mu \in \R \times T^*B$ is a Calder\'on projector for the normal operator $N(P)(\mu)$ of $P$  on $\Cinf(F_y;E)$.
\end{enumerate} 
\end{theorem}
Theorem \ref{thm:Calderon}(i) is  proven in Subsection \ref{ssec:Calderon}.
For statements (ii) and (iii)  on the $\phi$-symbol and the normal family of $C$, see Propositions \ref{prop:Calderon.symbol} and \ref{prop:Calderon.normal}. All pseudodifferential operators in this paper are classical, i.e.\ their symbols have complete expansions into positively homogeneous terms.

\subsubsection*{Choices} 
The different choices of $\calF$, i.e.\ of growth behavior at the singular boundary,  yield different boundary data spaces $\calB_{P,\calF}$. However, from an $L^2$ perspective these spaces are not very different: Recall the definition of $L^2_\phi$ and of the $\phi$-Sobolev space $H^k_\phi$ for a $\phi$-manifold (without \textBC-boundary), see Appendix \ref{sec:app mwc}. The space
$$  \calH= \bigoplus_{k=0}^{m-1}H^k_\phi(\dBCX;E) $$
is a natural space for boundary data in the $L^2$-setting, which for $m=1$ reduces to $L^2_\phi(\dBCX;E)$.
Any operator $C$ as constructed in Theorem \ref{thm:Calderon} is a bounded projection $\calH\to\calH$. We call its range the \textbf{$L^2$-boundary data space} of $P$ and denote it by $\calB_{P,L^2}$. The following corollary shows in particular that this space does not depend on the choice of $C$.
\begin{corollary} \label{cor:L2 closure of BC spaces}
Let $X$ and $P$ be as in Theorem \ref{thm:Calderon}.
If $\calF(\dBCX) \subset L^2_\phi(\dBCX)$ then $\calB_{P,\calF}\subset\calH$, and its closure is equal to $\calB_{P,L^2}$.
\end{corollary}
In particular, the closure of $\calB_{P,\calF}$ in $\calH$ is the same for all such choices of $\calF$. The smallest $\calF$ satisfying $\calF(\dBCX) \subset L^2_\phi(\dBCX)$ is $\Cdotsinf$, other choices include $\calF=x^\alpha\Cinf$ with $\alpha>\frac {\dim B}2+1$.
There is also a more general statement with $L^2_\phi$ and $\calH$ replaced by the weighted spaces $x^\beta L^2_\phi$, $x^\beta\calH$ for any $\beta\in\R$.

We now discuss questions of uniqueness and canonical choices of a Calder\'on projector. This is relevant for instance when considering parameter-dependent problems or problems invariant under a group action. In general, there is no canonical choice of Calder\'on projector $C$. This is reflected in the fact that the construction of $C$ in the proof of Theorem \ref{thm:Calderon} involves several choices. However, we have a uniqueness result:
\begin{proposition}\label{prop:Calderon.unique}
The construction used in the proof of
Theorem \ref{thm:Calderon} determines $C$ uniquely modulo $\Psi^{-\infty}_\phi(\dBCX;E^m)$.

The full $\phi$-symbol of $C$, i.e.\ the element $[C]$ of $\Psi^*_\phi(\dBCX;E^m)/\Psi^{-\infty}_\phi(\dBCX;E^m)$ fixed in this way, is determined constructively by the infinite order jet of $P$ at $\dBCX$, i.e.\ by its equivalence class modulo $\CdotBCinf(X)\Diff^m_\phi(X;E)$. Here $\CdotBCinf(X)$ is the space of smooth functions on $X$ vanishing to infinite order at $\dBCX$.
\end{proposition}
The normal family of $C$ cannot be expected to be uniquely determined since it is a Calder\'on projector for the normal family of $P$, and again there is no canonical choice of such a projector. So the space $\Psi^{-\infty}_\phi(\dBCX;E^m)$ cannot be replaced by $x\Psi^{-\infty}_\phi(\dBCX;E^m)$.
Proposition \ref{prop:Calderon.unique} is proven in Subsection \ref{ssec:Calderon.properties}.
Interestingly, the proof shows that it is much easier to construct the equivalence class $[C]$ than an operator $C$ which is an actual projection to the boundary data space.

With additional data there are also ways to distinguish special choices of $C$. For simplicity we consider the $L^2$-setting.
Since only the range of $C$ is determined by the given data $X$ and $P$ (and $\nu$), the  possible projections $C$ are parametrized by their kernels, which must be subspaces of $\calH$ complementary to the $L^2$-boundary data space of $P$.
There are (at least) two approaches in the literature how to choose additional data to fix a kernel and hence a projection, in the setting of a compact manifold with boundary:
\begin{enumerate}
 \item
 If an extension of $X$ to a closed manifold $\Xhat$, as well as extensions of the bundles $E,E'$ to $\Xhat$ and of $P$ to an elliptic, \emph{invertible} operator $\Phat$ on $\Xhat$ are given, then the boundary data space for $\Phat$ from the \lq other\rq\ side, $\Xhat\setminus X$, is complementary to the boundary data space of $P$ on $X$, and the corresponding projection is a  Calder\'on projector.
 \item
If one chooses additional geometric data at $\dX$ which defines an $L^2$ scalar product then one may look at the \emph{orthogonal} projector.\footnote{Such a construction seems quite unnatural except in the case of Dirac operators which involve choices of geometric data anyway.}
\end{enumerate}
The projectors in both cases are pseudodifferential.
The one in (i) is in the class described in Proposition \ref{prop:Calderon.unique}, but the one in (ii) will in general even have a different principal symbol.
We start by extending (i) to the singular setting. 
In general, it is not possible to extend $P$ to a $\phi$-elliptic \textit{differential} operator  on $\Xhat$ (see Section \ref{subsec:inv ext}), so we shall work with $\phi$-pseudodifferential operators. By the general theory,  the following conditions on a fully elliptic operator $\Phat\in\Psi^m_\phi(\Xhat)$ on a $\phi$-manifold $\Xhat$ are equivalent:\footnote{For the proof use that a fully elliptic $\Phat$ has a parametrix $Q\in\Psi^{-m}_\phi(\Xhat)$ so that the remainder terms  $Q\Phat-I$, $\Phat Q-I$ are in $x^\infty\Psi^{-\infty}_\phi(\Xhat)$ and are projections to the kernel and cokernel of $\Phat$.}
\begin{itemize}
 \item
 $\Phat$ is invertible in $\Psi^*_\phi(\Xhat)$, i.e.\ there is $Q\in\Psi^{-m}_\phi(\Xhat)$ so that $\Phat Q=Q\Phat = I$.
 \item
 $\Phat:H^m_\phi(\Xhat)\to L^2_\phi(\Xhat)$ is invertible.
 \item
 $\Phat: \calF(\Xhat)\to \calF(\Xhat)$ is invertible for any choice of admissible function space $\calF$.
\end{itemize}
If this is satisfied, we simply say that $\Phat$ is invertible. For the following theorem recall the definition of $\calH$ and $\calB_{P,L^2}$ before Corollary \ref{cor:L2 closure of BC spaces}.

\begin{theorem} \label{thm:calderon canonical inv extension}
Assume that $P$ has an extension to an invertible, fully elliptic $\phi$-$\Psi$DO $\Phat$ on a $\phi$-manifold (without \textBC-boundary) $\Xhat$ extending $X$ across $\dBCX$, acting between extended bundles $E,E'$.
Denote by  $\calB_{P,L^2}^+$, $\calB_{P,L^2}^-$  the
$L^2$-boundary data spaces of $P$ and of $\Phat_{|\Xhat\setminus X}$, respectively. Then
\begin{itemize}
 \item $\calB_{P,L^2}^+ \oplus\calB_{P,L^2}^- = \calH$.
 \item
The projection $C_{\Phat}$ in $\calH$ with range $\calB_{P,L^2}^+$ and kernel $\calB_{P,L^2}^-$ is a $\phi$-$\Psi$DO as  in Theorem \ref{thm:Calderon}, and its normal family is the Calder\'on projector analogously defined using the $\pm$-boundary data spaces of the normal family of $\Phat$.
\item
If $\calF$ is any admissible function space then $\calB_{P,\calF}^+ \oplus\calB_{P,\calF}^- = \calF(\dBCX;E)^m$ and $C_{\Phat}$ acts as a projection in $\calF(\dBCX;E)^m$
with range $\calB_{P,\calF}^+$ and kernel $\calB_{P,\calF}^-$.
\end{itemize}
\end{theorem}
Note that $\Phat$ is allowed to be a \emph{pseudo}differential operator.
See Section \ref{ssec:extensions} for a precise definition of `extension', which implies in particular that the restriction of $\Phat$ to $\Xhat\setminus X$ makes sense also in this case. Theorem \ref{thm:calderon canonical inv extension} is proven in Subsection \ref{ssec:Calderon}. For example, it applies to Laplace-Beltrami operators for metrics $x^{2c}g$ where $g$ is a $\phi$-metric, as in the examples (i)-(iii) above, see Subsection \ref{ssec:laplacian.extension}. See the paragraph around \eqref{eqn:basis phi-cotangent} for the definition of $\phi$-metrics.

We now consider (ii), i.e.\ orthogonal Calder\'on projectors.
For simplicity we only consider $m=1$, which includes the important class of Dirac type operators.
\begin{theorem} \label{thm:orthogonal}
 In the setting of Theorem \ref{thm:Calderon} assume a $\phi$-metric is given on $X$, as well as a hermitian metric on $E$. Assume  $m=1$.
Then the orthogonal projection in  $L^2_\phi(\dBCX;E)$ to $\calB_{P,L^2}$ is in $\Psi^0_\phi(\dBCX;E)$.
\end{theorem}
See Proposition \ref{thm:orth.Calderon} for a more precise statement and the proof, including the calculation of the principal symbol and normal family. Note that the principal symbol is, in general, different from that of the projector in Theorem \ref{thm:Calderon}.

\subsubsection*{Extensions and perspectives} Both the singular boundary $\dsX$ and the \textBC-boundary $\dBCX$ are allowed to have several connected components, see Definition \ref{def:bv-mfd}. For $\dsX$ this would correspond to several singularities. In fact, the dimensions of fibres and bases could vary from component to component, and
some components may have empty intersection with $\dBCX$, so that their fibres have empty boundary. These would correspond to interior singularities (in contrast to  boundary singularities). Assumption  \ref{UCNF} localizes to each component of $\dsX$, and at components not intersecting $\dBCX$ it
is equivalent to the operator $P$ being fully elliptic there, and is also a necessary condition for our theorem to be provable within the $\phi$-calculus.

On the other hand, it is not clear whether our theorems also hold without Assumption \ref{UCNF}. At least our methods do not carry over to that case, see Remark \ref{rem:need UCNF}.

A side result of our analysis is the following fact. Recall that $u\not\equiv0$ is called a shadow solution of $Pu=0$ if it is smooth and has zero boundary data at $\dBCX$.
\begin{proposition} \label{prop:shadow solutions}
 Let $X$ and $P$ be as in Theorem \ref{thm:Calderon}. Then any shadow solution  of $Pu=0$ is rapidly decreasing at the singular boundary $\dsX$; more precisely, it lies in the space $\CdotsBCinf(X;E)$ defined in Subsection \ref{ssec:function.spaces}.
\end{proposition}

Finally, we remark that, as in the classical constructions by Seeley and by H\"ormander, the operator $P$ only needs to be differential \emph{near} $\dBCX$, i.e.\ $P$ is an elliptic $\phi$-$\Psi$DO and there is a neighborhood $U$ of $\dBCX$ on which $P$ is a $\phi$-differential operator, and $\supp Pu \subset X\setminus U$ whenever $\supp u\subset X\setminus U$.

We have chosen to treat the boundary $\dBCX$ as in the classical literature on boundary value problems, while using the methods of geometric microlocal analysis at the singular boundary $\dsX$. An alternative approach would be to deal with $\dBCX$ using the zero-calculus, since if $P$ is an $m$th order differential operator which is elliptic uniformly to $\dBCX$ then $\rho^m P$ is an elliptic element of that calculus, where $\rho$ is a defining function for $\dBCX$. The zero calculus (for a manifold with boundary, i.e. $\dsX=\emptyset$) is treated in  \cite{MazMel:MERCSWACNC}, \cite{Maz:HCCCM} and \cite{Lau:PACCS}. 

\subsection{Outline of the paper}\label{ssec:outline}
In Section \ref{sec:setting} we introduce the setting in detail and explain how the examples mentioned above fit into this framework.
The proof of Theorem \ref{thm:Calderon} proceeds along the lines of Seeley's proof in the non-singular case, \cite{See66} and \cite{See69}, combined with the construction by H{\"o}rmander, \cite[Sec.\ 20]{Hor85}. In Section \ref{sec:outline} we give an outline of this construction. In Section \ref{sec:aug-ext-mod} we introduce the notions of augmenting, extending and modifying, leading to generalized extensions of spaces, bundles and operators.
These abstract notions help to clarify Seeley's construction, and we hope that they will be useful for extending our results  to different geometries. In  Section \ref{sec:inv.ext} we construct an invertible generalized extension $\Phat$ of $P$. Its properties imply Proposition \ref{prop:shadow solutions} as shown in Subsection \ref{ssec:pf prop shadow}. Then we use its inverse to construct a Calder\'on projector $C$  with the stated properties in Section \ref{sec:Calderon}. The construction of the Calder\'on projector from the inverse of $\Phat$ requires considering the transmission property for $\phi$-$\Psi$DOs. This is straightforward once we rephrase the transmission property in the language of conormal distributions (as opposed to $\Psi$DOs, as is usually done). This is carried out in Section \ref{ssec:transmission} and may be of independent interest. Theorems \ref{thm:Calderon} and \ref{thm:calderon canonical inv extension} and Corollary \ref{cor:L2 closure of BC spaces} are proved in Section \ref{ssec:Calderon}.
Theorem \ref{thm:orthogonal} follows from Theorem \ref{thm:Calderon} via standard properties of the $\phi$-$\Psi$DO calculus and by applying a formula for orthogonal projections as for instance in \cite{BLZ09}.

We deal with some parts of the construction in the greater generality of manifolds with \textBC-boundary in order to facilitate generalization to other types of singularities.

\subsubsection*{Acknowledgment} We thank the referee for his or her careful reading of the manuscript and helpful suggestions and comments.

\section{Setting and Examples} \label{sec:setting}

We now introduce the setting in detail. We use basic notions of manifolds with corners, which are collected in Appendix \ref{sec:app mwc}.
The setting is characterized by a class of spaces together with additional structure, which defines a class of partial differential operators on them, and a notion of ellipticity for these operators. All of these will be introduced in Subsection \ref{ssec:phi-geom}. Then, in Subsection \ref{ssec:function.spaces}, we introduce the type of function spaces we consider as well as the associated spaces of boundary data and shadow solutions.

\subsection{\texorpdfstring{$\phi$}{\straightphi}-{\scriptsize BC}-geometry}\label{ssec:phi-geom}
The geometric structure is described by the following definitions. It may be helpful for the reader to look at Figure \ref{fig:blow-up example} for illustration.

\begin{definition}\label{def:bv-mfd} A \textbf{manifold with} \textBC-\textbf{boundary}, or \textBC-\textbf{manifold}, is a smooth compact manifold with corners $X$, with a choice of a disjoint union of boundary hypersurfaces of $X$ designated as its \textBC-\textbf{boundary} and denoted by $\dBCX$. We then denote $\dsX = \overline{\partial X \setminus \dBCX}$, the union of the remaining boundary hypersurfaces, and call this the \textbf{singular boundary}.

A \textbf{$\phi$-manifold with} \textBC-\textbf{boundary}, or short $\phi$-\textBC-\textbf{manifold}, is a \textBC-manifold $X$ so that $\dsX=\bigcup_i \dsiX$ is also a disjoint union of boundary hypersurfaces of $X$, together with the following data.
\begin{enumerate}
 \item For each $i$ a fibration $ F_i - \dsiX \overset{\phi_i}\to B_i$ where the base $B_i$ is a smooth  closed manifold and the fibre $F_i$ is a smooth compact manifold with boundary, whose boundary corresponds to the \textBC-boundary $\dsBCiX := \dsiX \cap \dBCX$ of $\dsiX$:
\begin{equation}
\label{eqn:fibrations}
  \begin{tikzcd}[ampersand replacement=\&]
         F_i \arrow[r,dash] \& \dsiX\  \arrow[dr,"\phi_i"] \\[-4mm]
         \&\& B_i \\[-4mm]
         \partial F_i \arrow[uu,hook] \arrow[r,dash] \& \dsBCiX \arrow[uu,hook] \arrow[ur,swap,"\phiBCi" near start]
       \end{tikzcd}
\end{equation}
where $\phiBCi$ is the restriction of $\phi_i$ to $\dsBCiX$.

\item A boundary defining function $x$ for $\dsX$, i.e.\ $x:X\to\R_+=[0,\infty)$ is smooth, $\dsX=\{x=0\}$ and $dx\neq0$ at $\dsX$.
\end{enumerate}
\end{definition}
To simplify the notation we will assume for the most part that $\dsX$ is connected and then leave out the index $i$, and will point out adjustments for the disconnected case only when they are not obvious.
We will in general simply write $\phi$ instead of $\phiBC$ to keep the notation short.
The fibration $\phi$ of $\dsX$ can be extended to a product neighborhood of $\dsX$, and we will fix such an extension throughout.

Note that for a \textBC-manifold $X$, $\partial X = \dsX \cup \dBCX$ and any boundary hypersurface is contained in precisely one of $\dsX$ or $\dBCX$. The classical non-singular case (compact manifold with boundary) corresponds to $\dsX = \emptyset$. In the definition we allow the \textBC-boundary to be empty, in which case we obtain manifolds with corners resp. \textbf{$\phi$-manifolds} (in that case $\partial F=\emptyset$ also).
The class of $\phi$-manifolds was introduced by Mazzeo and Melrose in \cite{MM98}, and they are called manifolds with fibred boundary there.

For a $\phi$-\textBC-manifold $X$, only $\dsBCX$ can be a codimension two corner and there are no higher codimension corners. $\dsBCX$ is a closed manifold. If non-empty, the \textBC-boundary $\dBCX$ is a $\phi$-manifold whose boundary fibrations are given by the bottom line in \eqref{eqn:fibrations}.

The \textBC-boundary of a \textBC-manifold $X$ should be thought of as a boundary as it appears in classical boundary value problems, for example, while  $\dsX$
serves as a smooth model for the singularities. Different components $\dsiX$ of $\dsX$ may model different kinds of singularities, e.g.\ the fibres $F_i$ may have varying dimensions.
In some cases, e.g.\ in the example in Figure \ref{fig:blow-up example}, the analysis on $X$ models that on the singular space $\tilde X$ obtained by collapsing the fibres of $\phi$ to points, and conversely $X$ is obtained from a singular space $\tilde X$ by a suitable blow-up of its singular set.

Generally we consider objects on $X$, e.g.\ functions or metrics, which are defined on $X\setminus \dsX$ and smooth up to $\dBCX\setminus\dsBCX$, while they may be undefined at $\dsX$. One way to express smoothness up to the \textBC-boundary is via extensions across that boundary, so we define:

\begin{figure}
 \includegraphics[scale=1.1]{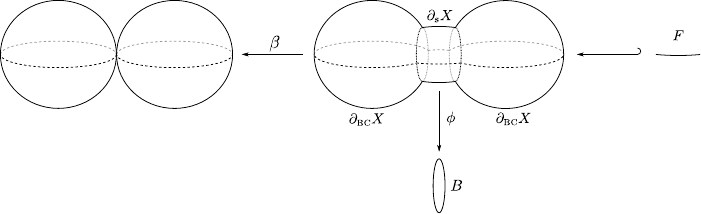}
 \caption{An example of a $\phi$-\textBC-manifold $X$ (exterior of the central picture) and how it arises via blow-up from a
 singular space (exterior of left picture, cf.\ Example \ref{ex:example touching domains}(i)). Here the boundary face $\dsX$ of $X$ created by the blow-up is fibred with base $B=S^1$ (bottom) and fibre $F=[-1,1]$ (right).}
 \label{fig:blow-up example}
\end{figure}

\begin{definition}\label{def:bv-ext} If $X$ is a manifold with non-empty \textBC-boundary, then a \textBC-\textbf{extension} of $X$ is a smooth compact manifold with corners $\Xhat$ containing $X$ and of the same dimension as $X$ and so that $\dBCX$ is an interior $p$-submanifold of $\Xhat$ and $\dsX \subset \partial \Xhat$.

If $X$ is a $\phi$-manifold with non-empty \textBC-boundary and with fibration $\phi: \dsX \rightarrow B$, then a $\phi$-\textBC-\textbf{extension} of $X$ is a \textBC-extension $\Xhat$ of $X$ so that $\Xhat$ is a $\phi$-manifold (without \textBC-boundary) whose  fibration $\phihat : \partial \Xhat \rightarrow B$ has the same base as $\phi$ and extends $\phi$, i.e., $\phi$ is the restriction of $\phihat$ to $\dsX$.
\end{definition}

Note that the fibres $F$ of $\phi$ are \textBC-manifolds (with $\partial_\BC F = F \cap \dBCX = \partial F$ and $\partial_\s F = \emptyset$) and that the fibres $\Fhat$ of $\phihat$ are \textBC-extensions of the fibres $F$ of $\phi$.
\textBC-- and $\phi$-\textBC-extensions may be obtained by doubling across the \textBC-boundary, see Lemma \ref{lem:mfd.extension}.

Whenever we express things in local coordinates, then near any point $q$ of $\dsX$ we will use \textbf{adapted coordinates} $x\geq0$, $y\in\R^b$, $z\in\R^f$ (or $z\in\R_+\times\R^{f-1}$ if $q \in \dsBCX$) centered at $q$, where $b=\dim B$, $f=\dim F$. Here $x$ is the given boundary defining function for $\dsX$, $y$ is pulled back from a local coordinate system for $B$, and $z$ are remaining coordinates, locally parametrizing the points in each fibre.

We now turn to differential operators on a $\phi$-\textBC-manifold $X$. This is an immediate extension of the definitions for $\phi$-manifolds from \cite{MM98} -- we simply take everything smooth up to the \textBC-boundary. The geometric data define the Lie algebra of $\phi$-\textbf{vector fields}:
\begin{equation}
\label{eqn:def phi-vector field}
\begin{aligned}
 \calV_\phi(X) = \{V\in \Cinf(X;TX):\, &Vx=O(x^2), \\  &V\text{ is tangential to the fibres of }\phi \text{ at }\dsX\}
\end{aligned}\end{equation}
They are unconstrained away from $\dsX$ but in adapted coordinates near a point of $\dsX$ are $\Cinf(X)$-linear combinations of the following vector fields:
\begin{equation}
 \label{eqn:basis phi-tangent}
 x^2\partial_x,\ x\partial_{y_i}, \partial_{z_j},\quad i=1,\dotsc,b,\ j=1,\dotsc,f
\end{equation}
Composing these vector fields and adding functions we obtain
$\phi$-\textbf{differential operators}: these have the form, near $\dsX$,
\begin{equation}
 \label{eqn:phi operator}
 P = \sum_{k+|\alpha|+|\beta|\leq m} a_{k,\alpha,\beta}(x,y,z) (x^2D_x)^k (xD_y)^\alpha D_z^\beta,\quad a_{k,\alpha,\beta} \text{ smooth, }
\end{equation}
where $D_x=\frac1i\partial_x$ etc., $\alpha,\beta$ are multi-indices and $m\in\N_0$ is the order of $P$, cf.\ \eqref{eqn:phi-op local coords}. The set of these operators is denoted $\Diff_\phi^m(X)$.  The class of operators acting between sections of vector bundles $E,E^\prime$ on $X$ is defined similarly and denoted $\Diff^m_\phi(X;E,E^\prime)$.
The $\phi$-\textbf{principal symbol} of $P\in\Diff^m_\phi(X)$ is the standard principal symbol in the interior, and near $\dsX$ where \eqref{eqn:phi operator} holds it is the function
$$\phisigma_m(P) =  \sum_{k+|\alpha|+|\beta| = m} a_{k,\alpha,\beta}(x,y,z) \tau^k \eta^\alpha \zeta^\beta$$
for $\tau\in\R$, $\eta\in\R^b$, $\zeta\in\R^f$.
Invariantly and globally, $\phisigma_m(P)$ can be made sense of as a function (or section of a suitable vector bundle of homomorphisms) on a rescaled cotangent bundle $\phiT^*X$ whose local basis near $\dsX$ is \eqref{eqn:basis phi-cotangent} below. $P$ is called $\phi$-\textbf{elliptic} if its $\phi$-principal symbol is invertible outside the zero section $(\tau,\eta,\zeta)=0$.

The \textbf{normal family} of $P\in \Diff^m_\phi(X)$ captures its behavior at the singular boundary. It is a family of differential operators
on the fibres $F_y$, parametrized by $y\in B$, $\tau\in\R$ and $\eta\in T_y^*B$, and having coefficients smooth up to $\partial F_y$. In adapted local coordinates where \eqref{eqn:phi operator} holds it is given by
\begin{equation}
\label{eqn:def normal op diff}
 N(P)(\tau;y,\eta) =  \sum_{k+|\alpha|+|\beta|\leq m} a_{k,\alpha,\beta}(0,y,z) \tau^k \eta^\alpha D_z^\beta
\end{equation}
where now $\eta\in\R^b$.\footnote{In the literature the notation $\hat N(P)$ is sometimes used for the normal family, with $N(P)$ denoting the normal \emph{operator}, where $\tau,\eta$ are replaced by differentiations $D_T$, $D_Y$ in auxiliary variables $T\in\R$, $Y\in\R^b$. We do not use the normal operator.}
If $\dBCX=\emptyset$ (e.g.\ for a \textBC-extension of our given $\phi$-\textBC-manifold) then $P$ is called \textbf{fully elliptic} if it is $\phi$-elliptic and $N(P)(\tau;y,\eta)$ is invertible for all $\tau,y,\eta$. Full ellipticity is equivalent to $P$ having a parametrix in the $\phi$-calculus, modulo errors that are smoothing and small near $\dsX$ in a suitable sense, and also to being Fredholm between $\phi$-Sobolev spaces. See Appendix \ref{sec:app phi ops} for details on the $\phi$-calculus and for an extension of the $\phi$-symbol and normal family to $\phi$-pseudodifferential operators.

\smallskip

The main motivation for considering $\phi$-differential operators is that geometric operators, e.g.\ the Laplacian and Dirac operators, for $\phi$-metrics are in this class (see \cite{MM98} and \cite{Mel90}). A $\phi$-\textbf{metric} is a  Riemannian metric on $X\setminus\dsX$ that, near $\dsX$, can be written as positive definite quadratic form in terms of
\begin{equation}
 \label{eqn:basis phi-cotangent}
\frac{dx}{x^2},\ \frac{dy_i}{x},\ dz_j,\quad i=1,\dotsc,b,\ j=1,\dotsc,f
\end{equation}
with smooth coefficients (where smoothness and positive definiteness hold up to $x=0$). More generally, the Laplacian for a metric of the form $x^{2c}g$, where $c\in\Z$ and $g$ is a $\phi$-metric, is of the form $x^{-2c}P$ where $P\in\Diff^2_\phi(X)$ is $\phi$-elliptic.

We will fix a $\phi$-\textBC-manifold $X$ throughout the main construction. To simplify the exposition we will also choose a background $\phi$-metric $g$ on $X$, smooth up to $\dBCX$. This is useful in two ways:
\begin{enumerate}
 \item The volume form $\dvol_g$ allows us to interpret Schwartz kernels of operators as distributions and avoid densities.
 \item  Having the $L^2$-space will allow us to talk about adjoints of operators, which will be useful in the construction.
\end{enumerate}
We also fix a $\phi$-vector field $\nu$ transversal to $\dBCX$ in order to define the boundary data map $\gamma$, see \eqref{eqn:first.gamma}. Given a $\phi$-\textBC-extension $\Xhat$ of $X$ we extend $\nu$ to a $\phi$-vector field on $\Xhat$. The flow of $\nu$ defines a trivialization
\begin{equation}
 \label{eqn:trivialize bd}
U \cong (-1,1)_\rho\times \dBCX
\end{equation}
(after scaling $\nu$ if needed)
of a neighborhood $U$ of $\dBCX$ in $\Xhat$ which is compatible with $\phi$ (that is, the image of a $\hat\phi$-fibre is a $\phi_\BC$-fibre times $(-1,1)$) and with the chosen boundary defining function $x$ and with $\rho\geq0$ in $X$, and $D_\nu=D_\rho$ then. In adapted coordinates we may and will take $z_1=\rho$.

\subsection{Function spaces and boundary data}\label{ssec:function.spaces}
Throughout, spaces of smooth functions vanishing at (parts of) the boundary of a \textBC-manifold $X$ will play an important role. We denote
 \begin{equation}\label{eqn:Cdotinf.BC}
  \CdotBCinf(X) = \big\{\, u \in \Cinf(X) \,:\, \text{$u$ vanishes to infinite order at $\dBCX$}\,\big\}
 \end{equation}
and similarly denote by $\Cdotsinf(X)$, $\CdotsBCinf(X)$ the spaces of smooth functions on $X$ vanishing to infinite order at $\dsX$ and at $\dsX \cup \dBCX$, respectively. If $x$ is a boundary defining function for $\dsX$, $\Cdotsinf(X) = x^\infty \Cinf(X)$ and $\CdotsBCinf(X) = x^\infty \CdotBCinf(X) = \Cdotinf(X)$, where $\Cdotinf(X)$ is the space of smooth functions vanishing to infinite order at $\partial X$. As all of these spaces are local $\Cinf(X)$-modules, we can define for instance $\CdotBCinf(X;E)$ for a vector bundle $E \rightarrow X$.

As mentioned in the introduction, we also need to specify growth conditions on functions on a \textBC-manifold $M$ near the singular boundary. These are encoded in function spaces which are admissible in the sense of the following definition. Recall from Appendix \ref{sec:app mwc} the definition of $\calA_s(M)$, the space of functions conormal at $\dsM$.

\begin{definition} \label{def:admissible function space} An assignment $\calF$ of a function space $\calF(M)$ to any \textBC-manifold $M$ is called \textbf{admissible} if it satisfies the following conditions, for every $M$:
\begin{enumerate}
 \item $\Cdotsinf(M) \subset \calF(M) \subset \calA_\s(M)$
 \item If $\dBCM\neq\emptyset$ and  $\Mhat$ is a \textBC-extension of $M$ then the restriction maps
        $$ \calF(\Mhat) \to \calF(M), \qquad \calF(M) \to \calF(\dBCM)$$
     are defined and surjective.
 \item $\calF(M)$ is a local $\Cinf(M)$-module.
\end{enumerate}
In the context of $\phi$-\textBC-manifolds we require (ii) to hold for $\phi$-\textBC-extensions and in addition:
\begin{enumerate}
 \item[(iv)]
 If $M$ is a $\phi$-manifold (i.e.\ $\dBCM=\emptyset$) then $\calF(M)$ is stable under $\Psi^*_\phi(M)$, i.e.\ $u\in\calF(M)$, $P\in\Psi^*_\phi(M)\Rightarrow Pu\in\calF(M)$.
\end{enumerate}

Given such an admissible choice $\calF$, we define in analogy to $\CdotBCinf(X)$:
 \begin{equation}\label{eqn:dfn.FdotBC}
  \calFdotBC(M) := \{u\in\calF(M):\,u \text{ vanishes to infinite order at }\dBCM\}
 \end{equation}
\end{definition}
Note that (i) implies that functions $u\in\calF(M)$ are smooth up to the \textBC-boundary, but unrestricted otherwise, away from $\dsM$. Also, functions in $\calF(M)$ can be paired with functions in $\Cdotsinf(M)$, with respect to a $\phi$-volume form. Condition (ii) is used in the construction of Calder\'on projectors.
Condition (iii) ensures that spaces of sections $\calF(M,E)$ of vector bundles $E\to M$ are defined, and they have analogous properties. Condition (iv) will be used for $M$ being either the \textBC-boundary or a $\phi$-\textBC-extension of a $\phi$-manifold with non-empty \textBC-boundary.

Admissible choices for $\calF$ include spaces of conormal functions $\calA^a$, polyhomogeneous functions $\calA^\calE$ or functions in $x^k \Cinf$, $k \in \Z$. Their definitions are recalled in Appendix \ref{sec:app mwc}. Here the conormality order $a$ and index family $\calE$ refer only to the singular boundary. Note that the order and index family do not change under restriction to the \textBC-boundary. This is different from the Sobolev order which would decrease by $\frac12$.

Let us check that the conditions on $\calF$ imply that $\gamma$ defined in \eqref{eqn:first.gamma} maps as stated: Let $u\in\calF(X)$. Since $\nu$ is a $\phi$-vector field, $D_\nu$ is a $\phi$-differential operator, so $D_\nu^k u\in\calF(X)$ by condition (iv) in Definition \ref{def:admissible function space}. Then by condition (ii) $D_\nu^ku_{|\dBCX}\in\calF(\dBCX)$, so $\gamma u \in \calF(\dBCX)^m$. Note that this clearly implies $\calB_{P,\calF} \subset \calF(\dBCX)^m$ as in \eqref{eqn:first.bdy.data}.

\subsection{Examples} \label{ssec:ex}
In this section we give some examples of settings to which our results apply. Examples (i)-(iii) elaborate on the examples mentioned in the introduction.
In all examples metrics of the form $x^{2c}g$ arise, with $c\in\{0,1,2\}$  and $g$ a $\phi$-metric, so the Laplacian is $x^{-2c}$ times an elliptic $\phi$-operator.

\subsubsection*{Example (i): Incomplete fibred cusps}\label{ex:example touching domains}
Consider the set
$X_0=\{(\xi,\eta)\in\R^m\times\R^k: |\xi|\leq |\eta|^2, |\eta|<1\}$ where $m,k\in\N$, whose boundary $ |\xi|= |\eta|^2$ has a singularity at $\xi=0,\eta=0$. We introduce quasihomogeneous polar coordinates by writing
$$ X = [0,1)\times \Sph^{k-1}\times \B^m\,,\ \
\beta: X\to X_0,\ (x,\omega,z) \mapsto (\xi=x^2z,\eta = x\omega)\,, $$
where $\Sph^{k-1}=\{\omega\in\R^k:\,|\omega|=1\}$ and $\B^m=\{z\in\R^m:\, |z|\leq 1\}$.\footnote{Alternatively, one may resolve $X_0$ by two standard blow-ups: first blow up the point $(\xi,\eta)=0$, then the intersection of the lift of the $\xi=0$ plane with the front face.}
Except for the non-compactness at $x=1$, the space $X$ is a $\phi$-\textBC-manifold, with $\dBCX=[0,1)\times \Sph^{k-1}\times \Sph^{m-1}$ corresponding to $|\xi|=|\eta|^2$ and $\dsX=\{0\}\times \Sph^{k-1}\times \B^m$, and $\phi:\dsX\to B=\Sph^{k-1}$ the projection, so the fibre is $\B^m$.
If $y$ are local coordinates on the sphere then the Euclidean metric $|d\xi|^2+|d\eta|^2$ pulls back to a positive definite quadratic form in  $dx$, $x dy$, $x^2 dz$, with coefficients smooth in $x\geq0$, hence is $x^4 g$ for a $\phi$-metric $g$ on $X$.

Any bounded domain in $\R^{m+1}$ whose boundary is smooth except for isolated outward-pointing cuspidal singularities is locally of this type with $k=1$.

If $m=1$ then $X_0$ is the complement of the two solid paraboloids $\{\pm\xi>|\eta|^2\}$. The same resolution and type of metric arises on $X_0=\R^{k+1}\setminus(\Omega_1\cup\Omega_2)$ where $\Omega_1,\Omega_2$ are domains with smooth boundary whose closures intersect in a single point $p$ and which are simply tangent there, see Figure \ref{fig:blow-up example} for an example with $k=2$.
Note that the geometry at infinity of $\R^{k+1}$ is also of $\phi$-type (see Example (iii)). However, the Laplacian is not fully elliptic at infinity, so if both $\Omega_i$ are bounded then $X$ does not satisfy Assumption \ref{UCNF} there. We expect our result to hold anyway. Our results do apply to the part of $X_0$ lying in a large ball, which has an additional \textBC-boundary.

We remark that there is an \lq obvious\rq\ extension of $X_0$ without \textBC-boundary, namely $\R^{m+k}$, and one may construct the
Calder\'on projector using an inverse of the Laplacian on $\R^{m+k}$ in the construction as outlined in Section \ref{sec:outline}.
However, it is not clear that $C$ is a $\phi$-$\Psi$DO then. (In \cite{Fri14}, it is shown that the boundary layer potentials which can be used to construct $C$ are elements of the full $\phi$-calculus only, but it is the small calculus we need here. \cite{Fri19} studies the same problem in the setting of half-spaces, thus in a different geometry.)
The point is that the complement $\R^{m+k}\setminus X_0$ has a very different geometric structure near the singularity. Our double of $X_0$ is really a 'thin' double, i.e.\ it still has a cusp (though without \textBC-boundary).

\subsubsection*{Example (ii): Domains in locally symmetric spaces}
Consider $X_0=\{w\in\C:\,|\Re w| \leq c,\ \Im w \geq 1\}$ for some $c>0$ with the hyperbolic metric $g_0=\frac{dw\,d\wbar}{(\Im w)^2}$. Introduce coordinates $x=\frac1{\Im w}$, $z=\Re w$ and compactify $X_0$ by adding a singular boundary $\dsX=\{x=0\}$ to obtain $X=[0,1]_x\times[-c,c]_z$.
The metric is $x^2\left(d(\frac1x)^2 + dz^2\right) = x^2g$ where $g=\frac{dx^2}{x^4} + dz^2$. This has the form of a $\phi$-metric without $y$-coordinates, so for the fibration whose base is a point. This makes $X$ into a $\phi$-manifold with $\dBCX=\{|z|=c\}\cup \{x=1\}$, apart from the corners at $\{1\}\times\{\pm c\}$ which could be rounded without affecting the geometry near $x=0$.
A similar geometry arises from domains in other locally symmetric spaces of $\Q$-rank one, if the boundary of the domain extends into the cusps.

\subsubsection*{Example (iii): Fibre bundles over cones}
Let $B\subset\Sph^{n-1}$ be a compact submanifold and $C=\{rw:\,w\in  B, r>1\}\subset\R^n$. In polar coordinates the metric on $C$ induced by the Euclidean metric is $g_C=dr^2+r^2g_0$, where $g_0$ is the induced metric on $B$.
To describe its behavior near infinity we introduce the coordinate $x=\frac1r$. Then the metric is
$\frac{dx^2}{x^4} + \frac{g_0}{x^2}$, which is a $\phi$-metric for the space $\Cbar$ defined as  $C$ with a copy of $B$ added  at $x=0$ as the singular boundary; thus, $\Cbar$ is the radial compactification of $C$. Since there are no $z$-coordinates we take $\phi$ to be the fibration with point fibres, so the base is $B$ and $\phi$ is the identity.

The product $X=\Cbar\times F$ with a compact Riemannian manifold $F$ with boundary
 is a $\phi$-\textBC-manifold, where $\phi:B\times F\to B$ is the projection and $\dBCX=\Cbar\times\partial F$, and the product metric is a $\phi$-metric. A simple example is the slab $\R^n\times[0,1]$.
 A closely related embedded example is obtained by thickening the conical set $C$ in the orthogonal direction, i.e.\ choosing  $C_\eps=\{p+z:\,p\in C, z\perp T_pC, |z|\leq\eps\}$ for  $\eps>0$. If $\eps$ is  sufficiently small then this is an $n$-dimensional submanifold with boundary of $\R^n$. It can be parametrized as $(1,\infty)\times N^\eps B\to C_\eps$,
$(r,(w,z)) \mapsto rw + z$ where $N^\eps B =\{(w,z)\in NB:\, |z|\leq\eps\}$ and $NB$ is the normal bundle of $B$ in $\Sph^{n-1}$.
Adding a copy of $N^\eps B$ at $x=0$ (i.e.\ $r=\infty$) as singular boundary again we obtain a $\phi$-manifold with $\phi:N^\eps B\to B$ the bundle projection with
fibre a ball, and it follows from standard calculations (see e.g.\ \cite{Gra:T}) that the Euclidean metric on $C_\eps$ is a $\phi$-metric.
An example of a compact space $X$ to which our results apply directly is obtained by \lq capping off\rq\ $\Cbar$ smoothly near the origin of $\R^n$.

Note that the closed half space, radially compactified, is not an example of a $\phi$-\textBC-manifold.
Its geometry at infinity is that of a $\phi$-manifold where the base is a half sphere, thus has boundary (and not the fibres, which are points).

\subsubsection*{Example (iv): Exterior problem for a smooth domain}
If $\Omega\subset\R^n$ is a bounded domain with smooth boundary then our results apply to the operator $-\Delta + 1$ on $\R^n\setminus\Omega$. Here the \textBC-boundary $\partial\Omega$ and the singular boundary (at infinity) are disjoint. The geometry at infinity is that of $C$ in Example (iii) with $B=\Sph^{n-1}$. Note that since the \textBC-boundary does not extend to infinity, Assumption \ref{UCNF} requires full ellipticity for the operator, which is satisfied by $-\Delta+1$, but not by $\Delta$.

\section{Outline of the Construction} \label{sec:outline}

For the convenience of the reader we give an outline of the classical construction of the Calder\'on projector, ignoring all technical details, in the case of a compact manifold with boundary (no singularities, no bundles), and then point to the various technical issues that arise both in the classical and in the singular setting. We carry out the construction of a Calder\'on projector as in \cite[20.1]{Hor85} but should also mention \cite[Ex.\ 1.9]{Gru96}, which is closer to the original construction in \cite{See66}.

Let $P$ be an elliptic differential operator of order $m$ on the compact manifold $X$ with boundary $\dX$. By doubling $X$ across its boundary we obtain a compact manifold without boundary,
$\Xhat \supset X$.

\subsection{Constructing a Calder\'on projector from an invertible extension}\label{ssec:outline Calderon constr}
Assume first that $P$ can be extended to an invertible elliptic differential operator $\Phat$ on $\Xhat$.
Similarly to \cite{Hor85}, but using the exact inverse of $\Phat$ rather than a parametrix, one can then explicitly construct a pseudodifferential Calder\'on projector for $P$ as follows. Identify a neighborhood of $\dX$ in $\Xhat$ with $(-1,1)_\rho\times\dX$, where $\rho=0$ corresponds to $\dX\subset\Xhat$. Consider the boundary data map (cf.\ \eqref{eqn:first.gamma})
 \begin{equation}\label{eqn:cauchy.data}\begin{array}{rlrcl}
  \gamma: & \Cinf(X) \to \Cinf(\dX)^m \,,
   &u &\mapsto &(u_{|\partial},D_\rho u_{|\partial},\dots, D_\rho^{m-1} u_{|\partial})
 \end{array}\end{equation}
where $m$ is still the order of $P$ and denote `$\delta$ extension from the boundary' by (ignoring density factors)
 \begin{equation}\label{eqn:delta.data}\begin{array}{rlrcl}
  \gamma^\star:& \Cinf(\dX)^m \to \calD'(\Xhat) \,,
   &U=(U_0,\dots,U_{m-1}) &\mapsto &\sum_{l=0}^{m-1} D^l_\rho\delta(\rho) \otimes U_l
 \end{array}\end{equation}
where $\delta(\rho)$ is the Dirac measure of $\dX$. For $u\in\Cinf(X)$, denote by $u^0$ its extension by $0$ to $\Xhat$. The construction now proceeds as follows:
 \begin{enumerate}
  \item The main point is to observe that $\Phat (u^0) = (Pu)^0 + \gamma^\star\kerningJp_P \gamma u$ for a differential operator $\Jp_P:  \Cinf(\dX)^m\to \Cinf(\dX)^m$ of order $m-1$ determined by $P$, since extending $u$ by zero introduces a jump singularity (hence the letter $\Jp$) which yields (normal derivatives of) delta distributions when differentiated. If $Pu=0$ then the first term on the right vanishes.
  \item If $Pu = 0$, applying $\Phat^{-1}$ and taking boundary data yields $\gamma u = \gamma\Phat^{-1}\gamma^\star\kerningJp_P \gamma u$. Here $\gamma$ is understood as the limit when approaching $\dX$ from the interior of $X$. Denoting by $\calB$ the boundary data space defined in \eqref{eqn:cauchy data space} we see that the operator
         \begin{equation}\label{eqn:Calderon formula intro}
          C := \gamma\Phat^{-1}\gamma^\star\kerningJp_P : \ \calB\to\calB
         \end{equation}
      is the identity.
      Now the idea is to show that the formula defining $C$ defines an operator $\Cinf(\dX)^m\to \Cinf(\dX)^m$. Note that this would be given as the composition
        \begin{equation}\label{eqn:Calderon formula intro 2}
         \begin{array}{rccccccccl}
         \qquad C: &\Cinf(\dX)^m& \stackrel{\,\Jp_P}\to & \Cinf(\dX)^m & \stackrel{\gamma^\star}\to &\calD'(\Xhat)  &  \stackrel{\Phat^{-1}}\to & \calD'(\Xhat) & \stackrel{\gamma}\to & \Cinf(\dX)^m \\[1mm]
         & U &\mapsto& \Jp_PU  &\mapsto& f & \mapsto & v &\mapsto & \gamma v=CU
         \end{array}\end{equation}
      where all maps are well-defined except the last, since restriction of general distributions does not make sense.
      However, $\gamma$ only needs to be applied to $v=\Phat^{-1}f$ where $f=\gamma^\star\kerningJp_P U$ and $U\in\Cinf(\dX)^m$.
      Now by standard $\Psi$DO theory, the operator $\Phat^{-1}$ is pseudodifferential, hence pseudolocal. Since $f$ is smooth (in fact, zero) in $\interior{X}$, so is $v$, and although $v$ is singular at $\dX$, a further argument  involving the transmission property of $\Phat^{-1}$ shows that $v_{|\interior{X}}$ extends to a smooth function on $X$, so in this sense $\gamma v$ is defined, and that the operator $C:U\mapsto \gamma v$ is pseudodifferential.
      Finally, it remains to check that $CU\in\calB$ for all $U\in\Cinf(\dX)^m$.
      This is clear since $\Phat v=f$ and $f=0$ in $\interior{X}$, so $Pv=0$ in $\interior{X}$. Since $v$ extends smoothly to $\dX$ it follows that $\gamma v\in\calB$.
  \item The orthogonal Calder\'on projector $C_o$ defined by additional geometric data is given by the formula $C_o = C (\Id + C - C^\star)^{-1}$, where the adjoint $C^\star$ also depends on the geometric data;
      see \cite[Lem.\ 3.5]{BLZ09} for instance. It follows from standard $\Psi$DO theory again that $C_o$ is pseudodifferential as well.
 \end{enumerate}

\subsection{Constructing an invertible extension}\label{subsec:inv ext}
Above we assumed the existence of an invertible extension $\Phat$. However, such an invertible extension might not exist, at least not as a differential operator.

One obstruction is that the equation $Pu=0$ might have shadow solutions as
in Proposition \ref{prop:shadow solutions}.
Then by ellipticity $u$ is in $\dot C^\infty(X)$, the space of smooth functions vanishing to infinite order at $\dX$, and any differential operator $\Phat$ extending $P$ would have $u^0$ in its kernel.
Another obstruction to the existence of $\Phat$ is topological: An elliptic differential operator $P$ need not  have an {elliptic} extension to $\Xhat$ as differential operator.
\footnote{For example, consider the operator $P=\frac12(\partial_x+i\partial_y)$ on the unit disk $X\subset\R^2$. The double of $X$ is the 2-sphere $\Sph^2$, and there is no scalar elliptic first order scalar differential operator on $\Sph^2$  since its principal symbol  would furnish a trivialization $T^*\Sph^2\to\C\setminus\{0\}$ of the cotangent bundle of $S^2$, which does not exist by the `hairy ball theorem'. Of course this obstruction can be overcome by extending the trivial bundle on $X$ as the antiholomorphic 1-form bundle on $\Sph^2=\C P^1$, then the $\overline\partial$ operator taking values in sections of this bundle is an elliptic extension.}

Both of these obstacles can be overcome by constructing $\Phat$ as an elliptic \emph{pseudo}\-differen\-tial operator and by generalizing the notion of extension to also allow enlarging the bundles  and adding certain smoothing operators. Here is an outline of such a construction close to that of Seeley (\cite{See66} and the appendix of \cite{See69}) and used in this paper:
 \begin{enumerate}
 \item[(i')] Setting $\Pbar =  
 \begin{pmatrix} 0&P^\star\\P&0\end{pmatrix}$ one obtains a formally self-adjoint elliptic differential operator on $X$.
 We call this an augmentation of $P$; here one needs to choose an auxiliary background metric on $X$. Selfadjointness allows to connect the symbol of $\Pbar$ at $\dX$ via a homotopy through elliptic symbols to $|\xi|^m$ times the identity, and this  can be used to define an elliptic \emph{pseudo}differential operator $\Phat_0$ on $\Xhat$ that extends $\Pbar$. Augmenting again one obtains a self-adjoint elliptic operator $\Ptilde$ on $\Xhat$.
 \item[(ii')]
  If $\Ptilde$ is invertible then the construction above yields a Calder\'on projector for the restriction of $\Ptilde$ to $X$, and this easily yields one for $P$.
  \item[(iii')]
If $\Ptilde$ is not invertible the strategy is to perturb $\Ptilde$ to make it invertible, without changing the boundary data space $\calB_\Ptilde$.
One idea to make $\Ptilde$ invertible is to add an orthogonal projection to its kernel; however, this may change $\calB_\Ptilde$. Instead one can add an orthogonal projection $\Picomp$ to any subspace $\Vcomp$ complementary to the range of $\Ptilde$.
If $\Ptilde$ has no shadow solutions then $\Vcomp$ can be chosen to consist of functions supported in $X^-=\overline{\Xhat\setminus X}$, and then $\calB_{\Ptilde+\Picomp}=\calB_{\Ptilde}$, so we are done.

\item[(iv')]
 In general, if $\Vsh\subset\ker \Ptilde$ is the space of shadow solutions and $\Pish$ is the orthogonal projection to $\Vsh$ then $\Ptilde+\Pish$ has no shadow solutions and $\calB_{\Ptilde+\Pish}=\calB_{\Ptilde}$. So choosing $\Vcomp$ as in (ii'), but for $\Ptilde+\Pish$ instead of $\Ptilde$, we get that $\Phat':=\Ptilde+\Pish+\Picomp$ is invertible and has the same boundary data space as $\Ptilde$.
\end{enumerate}
Note also that the boundary data space is only defined for operators $\Xhat$ that may be restricted to $X$. This restriction is defined for $\Pish$ and $\Picomp$ since $\Vsh$, $\Vcomp$ consist of functions supported on one side of $\dX$, but not for the projection to the kernel.

The details about this generalized extension procedure are given in Section \ref{sec:aug-ext-mod}.

\subsection{Additional issues in the singular setting}
In our singular setting everything needs to be made to work within the $\phi$-$\Psi$DO calculus.
First, doubling $X$ (and the fibration $\phi$) across its \textBC-boundary we obtain a $\phi$-manifold $\Xhat$ without \textBC-boundary.
Once we have constructed an invertible generalized extension $\Phat'$ on $\Xhat$ which is a fully elliptic $\phi$-$\Psi$DO, we can apply the standard $\phi$-calculus of Mazzeo and Melrose \cite{MM98} to prove that its inverse is also a $\phi$-$\Psi$DO. After generalizing the transmission condition to the $\phi$-setting we can carry out steps (i)-(iii) essentially as before.

The main new ingredient in our singular setting is the following: In order to ensure that the inverse of $\Phat'$ is a $\phi$-$\Psi$DO, we need to construct $\Phat'$ in such a way that it is \textbf{fully} elliptic, i.e.\ that its normal family is invertible. For this one would like to do steps (iii'), (iv') above for each fibre and each operator in the normal family separately. However, the kernels of the normal family will in general not define a vector bundle.
This problem can be overcome for step (iii') since there is some flexibility here, but not for step (iv'), and it is at this point where we need to make the additional Assumption \ref{UCNF} about the  unique continuation property of the normal family.

\section{Augmenting, Extending and Modifying Operators} \label{sec:aug-ext-mod}

As explained in Section \ref{sec:outline} a central step in the construction of a Calder\'on projector for $P\in\Diff^m_\phi(X;E,E')$ is extending $P$ to an operator on a \textBC-extension of $X$, e.g.\ its \textBC-double. However, such an extension might not exist (cf.\ Subsection \ref{subsec:inv ext}), so we need to generalize the notion of extension. This generalization has three ingredients:
\begin{itemize}
 \item
 \textbf{Augmenting} the bundles $E,E'$ to larger bundles, and correspondingly the operator.
 \item Adding a projection onto a finite-dimensional subspace of $\CdotBCinf(X;E)$ (respectively $\calFdotBC(X;E)$) to an operator. We call this a \textbf{modification} of the operator.
 \item
 \textbf{Extending} the space $X$ to a $\phi$-manifold without \textBC-boundary, and correspondingly bundles and operators.
\end{itemize}
Augmentation is used to create formally self-adjoint operators, which then can be extended. Modification is needed to deal with operators having shadow solutions, since these cannot have an invertible extension.
See Section \ref{sec:inv.ext} for more details.

Combinations of these also occur; a combination of all three is called a \textbf{generalized extension}.

In this section we explain these notions in detail and show how to obtain a Calder\'on projector for $P$ from a Calder\'on projector for a modified augmentation of $P$.

Most of the arguments in this section are of a general functional analytic nature, about linear operators on some function space. With the application to $m$-th order differential operators in mind we  fix $m\in\N$, a boundary data map $\gamma$ and an admissible choice $\calF$ of function space, see Subsection \ref{ssec:function.spaces}. Thus, we also have the $\calF$-boundary data space $\calB_{T,\calF}$ of order $m$ of a linear operator $T : \calF(X;E) \to \calF(X;E')$ and the notion of an $\calF$-Calder\'on projector for $T$.

We adopt the following conventions for notation: a bar is used for an augmented object, e.g.\ $\Ebar$, and a hat is used for an extended object, e.g.\ $\Ehat$ (and also for an augmented extended object).

\subsection{Augmentations} \label{ssec:augmentations}
In this subsection $X$ denotes a manifold with \textBC-boundary, see Definition \ref{def:bv-mfd}. No $\phi$-structure is needed here.

If $E,F\to X$ are vector bundles then we can form the direct (or Whitney) sum $\Ebar=E\oplus F\to X$. It will be useful to phrase the relation of $E$ and $\Ebar$ without mentioning $F$, using the inclusion $E\hookrightarrow \Ebar$ and projection $\Ebar\twoheadrightarrow E$ with kernel $F$ instead:
\begin{definition}[Augmentations] \label{def:augmentation}
\mbox{}
\begin{enumerate}
 \item[(a)] Let $E\to X$ be a vector bundle. An \textbf{augmentation} of $E$ is a vector bundle $\Ebar\to X$ together with vector bundle maps
 $$ \iota: E \hookrightarrow \Ebar,\quad \pi: \Ebar \twoheadrightarrow E,\quad
\pi \circ \iota = \id_E $$
where $\iota$ is injective and $\pi$ is surjective.
 \item[(b)] Let $E,E'\to X$ be vector bundles and $T:\calF(X;E) \to \calF(X;E')$ be a linear operator.
 An \textbf{augmentation} of $T$ is a linear operator $\Tbar : \calF(X;\Ebar) \to \calF(X;\Ebar')$ between augmentations $\Ebar, \Ebar'$ of $E, E'$ for which
\begin{equation}
 \label{eqn:def augment op}
 \pi' \Tbar = T\pi\,,\quad \Tbar \iota = \iota' T \,.
\end{equation}
Here $\iota$ also denotes the induced map $\calF(X;E) \to \calF(X;\Ebar)$, and similarly for $\pi, \iota', \pi'$.
\end{enumerate}
\end{definition}
If we identify $E$ with $\iota E$ then $\Ebar=E\oplus F$ where $F=\ker\pi$, and similarly $\Ebar'=E'\oplus F'$. Then an
augmentation of $T$ is a direct sum (or block diagonal) operator
\begin{equation}
 \label{eqn:augment op direct sum}
\Tbar = T \oplus  U: E \oplus F \to E' \oplus F'
\end{equation}
 for some operator $ U:F\to F'$.
 We will use this characterization in this section.
 Later it will be useful to have both descriptions available.

We now consider the effect of augmentation on boundary data spaces and Calder\'on projectors.

\begin{proposition}
 \label{prop:augment B C}
  Let $\Tbar : \calF(X;\Ebar) \to \calF(X;\Ebar')$ be an augmentation of the linear operator $T : \calF(X;E) \to \calF(X;E')$. Then the maps $\iota: E \hookrightarrow \Ebar$, $\pi: \Ebar \twoheadrightarrow E$ induce an injection and surjection
 $$  \iota: \calB_{T,\calF} \hookrightarrow \calB_{\Tbar,\calF},\quad \pi: \calB_{\Tbar,\calF} \twoheadrightarrow \calB_{T,\calF} \,. $$
 Also, if $\Cbar: \calF(\dX;\Ebar)^m\to \calB_{\Tbar,\calF}$ is a Calder\'on projector for $\Tbar$ then
 $$ C = \pi \Cbar \iota: \calF(\dX;E)^m\to \calB_{T,\calF} $$
 is a Calder\'on projector for $T$.
\end{proposition}
\begin{proof}
 If we represent $\Ebar=E\oplus F$ and $\Tbar=T\oplus  U$ as in \eqref{eqn:augment op direct sum} then clearly $\calB_{\Tbar,\calF} = \calB_{T,\calF} \oplus \calB_{U,\calF}$, which implies the first claim. Also $\calF(\dX;\Ebar) = \calF(\dX;E)\oplus\calF(\dX;F)$, so the
 $\calF(\dX;E)^m\to \calB_{T,\calF}$ part of $\Cbar: \calF(\dX;E)^m\oplus\calF(\dX;F)^m\to \calB_{T,\calF}\oplus\calB_{U,\calF}$, which is $\pi \Cbar \iota$, is a projection to $\calB_{T,\calF}$.
\end{proof}

\subsection{Modifications} \label{ssec:modifications}
As in Subsection \ref{ssec:augmentations}, in this subsection $X$ denotes a manifold with \textBC-boundary.
Let $E\to X$ be a vector bundle. Operators in this subsection map sections of $E$ to sections of $E$ (so $E'=E$ in the previous notation).
As before we fix $m\in\N$, a boundary data map $\gamma$ of order $m$ and an admissible $\calF$.

\begin{definition}[Modifications]\label{def:modification}
 Let $T : \calF(X;E) \to \calF(X;E)$ be linear. A \textbf{modification} of $T$ is an operator $T+\Pi$ where
 $\Pi$ is a projection in $\calF(X;E)$ with
\begin{equation}
 \label{eqn:modification}
 \rg \Pi \subset  \ker T \cap \ker\gamma \,.
\end{equation}
\end{definition}

\begin{proposition}\label{prop:modification B}
 Let $T$ be as in Definition \ref{def:modification} and $T+\Pi$ be a modification. Assume moreover that  $\rg\Pi \cap \rg T = \{0\}$. Then
 $$ \calB_{T+\Pi,\calF} = \calB_{T,\calF}\,. $$
If, in addition, equality holds in \eqref{eqn:modification} then
\begin{equation}
 \label{eqn:modification no kernel}
 \ker (T+\Pi) \cap \ker\gamma = \{0\} \,.
\end{equation}
The condition $\rg\Pi \cap \rg T = \{0\}$ is satisfied, for example, if
$T$ restricts to an operator $\Cdotsinf(X;E) \to \Cdotsinf(X;E)$ and satisfies $\ker T \cap \ker \gamma \subset \Cdotsinf(X;E)$ and
 \begin{equation} \label{eqn:form selfadjoint}
 \langle Tu,v\rangle = \langle u,Tv \rangle\quad \text{ for }  u\in \Cdotsinf(X;E), v\in\calF(X;E), \gamma u=0 \,.
\end{equation}
\end{proposition}

\begin{remark}\label{rem:shadow solutions}
If $P$ is an elliptic differential operator of order $m$ then $\ker P \cap \ker \gamma = \ker P \cap \calFdotBC(X;E)$ by ellipticity.
The proposition says that these shadow solutions can be removed by adding to $P$ a projection $\Pi_{sh}$ onto $\ker P \cap \ker \gamma$ and that this does not alter the boundary data space, assuming the additional condition $\rg \Pi_{sh}\cap\rg P=\{0\}$, which follows from formal self-adjointness.
\end{remark}
The idea is that a `shadow solution' $u \in \rg \Pi$ is brought to light in $T+\Pi$ since $(T+\Pi)u=\Pi u=u$. The additional condition is needed to ensure that no new shadow solutions (or other solutions of $(T+\Pi)u=0$ altering $\calB$) are created.

\begin{proof}
Proof of $\calB_{T,\calF}\subset\calB_{T+\Pi, \calF}$:
 If $U\in\calB_{T,\calF}$ then $U=\gamma u$ with $Tu=0$. Let $v=u-\Pi u$.
 We claim that $(T+\Pi)v=0$ and $\gamma v=\gamma u$, which implies $U\in\calB_{T+\Pi, \calF}$. First, $\Pi v=0$ implies $(T+\Pi)v=Tv=Tu-T\Pi u=0$ since $Tu=0$ and $T_{|\rg\Pi}=0$. Second, $\gamma_{|\rg\Pi}=0$ implies $\gamma v=\gamma u$.

Proof of $\calB_{T+\Pi,\calF}\subset\calB_{T,\calF}$: If $U\in\calB_{T+\Pi,\calF}$ then $U=\gamma u$ with $(T+\Pi)u=0$. Then $Tu=-\Pi u$, and $\rg \Pi\cap\rg T=\{0\}$ implies $Tu=0$, so $U\in\calB_{T,\calF}$.

Proof of $\ker (T+\Pi) \cap \ker\gamma = \{0\}$: If $(T+\Pi)u=0$ and $\gamma u=0$ then $Tu=0$, $\Pi u=0$ as before, so $u\in \ker T\cap\ker\gamma=\rg\Pi$. Since also $u\in\ker\Pi$ and $\Pi$ is a projection it follows that $u=0$.

Proof of $\rg \Pi \cap \rg T = \{0\}$: Suppose that $T$ restricts to an operator $\Cdotsinf(X;E) \to \Cdotsinf(X;E)$ and let $u \in \Cdotsinf(X;E)$, $v \in \calF(X;E)$. Then also $Tu \in \Cdotsinf(X;E)$ and $Tv \in \calF(X;E)$ and because the $L^2$-scalar product extends to a pairing of $\Cdotsinf(X;E)$ and $\calF(X;E)$, the scalar products $\langle Tu, v \rangle$ and $\langle u, Tv \rangle$ are well-defined.
Now additionally assume $\ker T \cap \ker \gamma \subset \Cdotsinf(X;E)$ and \eqref{eqn:form selfadjoint}.
If $u\in\rg\Pi\cap\rg T$ then $u\in \rg\Pi\subset\ker T\cap \ker\gamma\subset \Cdotsinf(X;E)$ and $u=Tv$ for some $v\in\calF(X;E)$, so $Tu=0$ gives $0=\langle Tu,v\rangle = \langle u,Tv \rangle = \|u\|^2$, so $u=0$.
\end{proof}

\subsection{Extensions}\label{ssec:extensions}
\subsubsection{Extensions of spaces and bundles} \label{sssec:extspaces}

Given a \textBC-extension $\Xhat$ of $X$ as  in Definition \ref{def:bv-ext}
 we will sometimes denote
$$ X^+ = X\,,\quad X^- = \overline{\Xhat\setminus X}\,.$$
We need to show that ($\phi$-\textBC-)extensions exist.
One way to prove this is by doubling across the \textBC-boundary. For a compact manifold with non-empty boundary, $M$, this is standard: join two copies $M^\pm$ of $M$ along their identical boundaries: $\hat{M} = M^- \sqcup_{\partial M} M^+$ and use a collar neighborhood of $\partial M$ to define the smooth structure of $\Mhat$ near $\partial M$, yielding an extension of $M$. This extends to \textBC-manifolds, and also to $\phi$-\textBC-manifolds:

\begin{lemma}\label{lem:mfd.extension} Let $X$ be a $\phi$-manifold with non-empty \textBC-boundary and let $\Xhat = X^- \sqcup_{\dBCX} X^+$ be its double across $\dBCX$.
Then the $\phi$-structure of $X$ extends canonically to a $\phi$-structure on $\Xhat$ by doubling the fibres, making $\Xhat$ a $\phi$-\textBC-extension of $X$ in the sense of Definition \ref{def:bv-ext}.
\end{lemma}

\begin{proof}
We first choose a collar neighborhood of the \textBC-boundary compatible with $\phi$, that is, a diffeomorphism from a neighborhood $U\subset X$ of $\dBCX$
\begin{equation}
 \label{eqn:collar nbhd1}
 S:U\to[0,1)\times\dBCX;\quad S_{|\dBCX} = 0\times\id \,,
\end{equation}
which at $\dsX$ respects the fibres of $\phi$ (this makes sense since a diffeomorphism \eqref{eqn:collar nbhd1} must map $U\cap\dsX$ to $\dsBCX\times[0,1)$), as in \eqref{eqn:trivialize bd} but only for $\rho\geq0$.

The smooth structure on $\Xhat$ is inherited from that of $X$ away from $\dBCX$, and near $\dBCX$ is defined by choosing a collar neighborhood
$S$ as above and demanding the map
 \[ \widehat S:  \widehat U = U^- \sqcup_{\dBCX} U^+ \longrightarrow  (-1,1)\times\dBCX \]
to be a diffeomorphism, where $\widehat S$ is defined on $U^+$ as $S$ and on $U^-$ as $S$ followed by sign reversal in the first coordinate.
Then, defining
 \[ \phihat : \dsX^- \sqcup_{\dsBCX} \dsX^+ = \partial \Xhat \longrightarrow B \]
to be $\phi$ on both $\dsX^\pm$, we see that on $\widehat U$ it maps $\widehat S^{-1}(\rho,p)\mapsto \phi(p)$ for all $(\rho,p)\in(-1,1)\times\dsBCX$, hence is still a fibration. By a similar reasoning the boundary defining function $x$ for $\dsX\subset X$ on each copy $X^\pm$ yields a boundary defining function for $\partial\Xhat\subset \Xhat$.
\end{proof}

We emphasize that at $\dsX$ the doubling `happens in the fibres' only. The base $B$ is a manifold without boundary.

\begin{definition}
Let $E$ be a vector bundle over a $\phi$-manifold with non-empty \textBC-boundary, $X$, and let $\Xhat$ be a \textBC-extension of $X$.
An \textbf{extension} of $E$ is a vector bundle $\Ehat$ over $\Xhat$ so that $\Ehat_{|X} = E$.
\end{definition}

Note that if we take for $\Xhat$ the double of $X$ across $\dBCX$ then any vector bundle $E$ has an extension to $\Xhat$: simply use $E$ on $X^-$ and glue using
a trivialization of $E$ on a collar neighborhood of $\dBCX$. Nevertheless, we will generally use any \textBC-extension $\Xhat$ and then assume that the bundles under consideration extend as well.

\subsubsection{Extensions of operators}
If $\Xhat$, $\Ehat$ are (\textBC-)extensions of $X$, $E$, then any section $\uhat\in\calF(\Xhat;\Ehat)$ can be restricted to a section $\uhat_{|X}\in\calF(X;E)$, by the assumptions on $\calF$. If $u=\uhat_{|X}$ then we also say that $\uhat$ extends $u$. For operators on $\Xhat$ to be restrictable to $X$ is an extra condition.

\begin{definition}\label{def:extension operator}
Let $\Xhat$ be a \textBC-extension of $X$ and $\Ehat,\Ehat^\prime\to\Xhat$ be extensions of vector bundles $E,E^\prime\to X$.
We say that a linear operator $\hat{T} : \calF(\Xhat;\Ehat) \to \calF(\Xhat;\Ehat^\prime)$
\textbf{restricts} to $X$ if for all $\uhat\in\calF(\Xhat;\Ehat)$
\begin{gather}
 \label{eqn:extension of ops 1}
 (\That\uhat)_{|X} \text{ only depends on } \uhat_{|X}\\
 \label{eqn:extension of ops 2}
 \supp \uhat \subset X \Rightarrow \supp \That \uhat \subset X\,.
\end{gather}
In this case we denote the restriction by
$$\That_X: \calF(X;E) \to \calF(X;E^\prime)\,.$$
Given $T:  \calF(X;E) \to \calF(X;E^\prime)$ we say that $\That$ \textbf{extends} $T$ if $\That_X=T$.
\end{definition}

Note that $\That_X$ is unique if it is defined since restriction to $X$ is surjective $\calF(\Xhat;E)\to\calF(X;E)$.
Condition \eqref{eqn:extension of ops 1} says that for $u\in\calF(X;E)$,  $\That_X u :=(\That \uhat)_{|X}$
is independent of the choice of extension  $\uhat$ of $u$, so $\That_X$ is well-defined.
This is clearly equivalent to condition \eqref{eqn:extension of ops 2} with $X$ replaced by $X^-=\overline{\Xhat\setminus X}$. Therefore, $T$ restricts to $X$ if and only if it restricts to $X^-$.\footnote{In order to define $\That_X$, only condition \eqref{eqn:extension of ops 1} would be needed. However, we add condition \eqref{eqn:extension of ops 2} for the sake of symmetry, since it will be useful when considering adjoints.}
Also, note that \eqref{eqn:extension of ops 2}  is equivalent to
$$ \That_X : \calFdotBC(X;E) \to \calFdotBC(X;E^\p) \,, $$
with $\calFdotBC(X;E)$ as in \eqref{eqn:dfn.FdotBC}, since we can identify, via extension by zero,
\begin{equation} \label{eqn:CdotBCinf identific}
\calFdotBC(X;E)
= \{\uhat\in \calF(\Xhat;\Ehat):\, \supp \uhat\subset X \} \,.
\end{equation}
We close with additional remarks regarding Definition \ref{def:extension operator}:
\begin{remark}\label{rmk:extension operator}
 \mbox{}
 \begin{enumerate}
  \item If $\That$ is given in terms of a Schwartz kernel $K$, a distribution on $\Xhat\times\Xhat$ (as for $\Psi$DOs), then $\That$ restricts to $X$ if and only if $\supp K \subset X^2 \cup (X^-)^2$, i.e if and only if $K$ has block-diagonal structure with respect to the $\pm$ decomposition of $\Xhat$.
  So in this case (which is all we care about) Definition \ref{def:extension operator} is independent of the choice of function space $\calF$.
  \item If $\That$ restricts to $X$ then so does its adjoint  (with respect to any smooth measure on $\Xhat$ and hermitian metrics on the bundles).
  \item A projection $\Pi$ to a finite dimensional subspace $K\subset\calF(\Xhat;\Ehat)$ restricts to $X$ if
      \begin{itemize}
       \item either all $u\in K$ have support in $X$, or all $u\in K$ have support in $X^-$, and
       \item $\Pi$ is an orthogonal projection with respect to some metrics on $\Xhat$ and $\Ehat$.
      \end{itemize}
      One way to see this is that the Schwartz kernel of $\Pi$ is $\sum_j u_j\otimes \overline{u_j}$, for an orthonormal basis $(u_j)_j$ of $K$.
 \end{enumerate}
\end{remark}

\section{Constructing an Invertible Generalized Extension}\label{sec:inv.ext}

In this section we prove the following theorem.

\begin{theorem} \label{thm:construction inv ext}
 Let $X$ be a $\phi$-\textBC-manifold with non-empty \textBC-boundary and suppose $P\in\Diff_\phi^m(X;E,E')$ is a $\phi$-elliptic differential operator on $X$ satisfying Assumption \ref{UCNF} stated below. Let $\Xhat$ be a $\phi$-\textBC-extension of $X$ and assume that the bundles $E$ and $E^\p$ extend to $\Xhat$.
Then there is an augmented extension $\Ehat\to\Xhat$ of the bundles $E,E'\to X$ and an operator $\Phat\in\Psi^m_\phi(\Xhat;\Ehat)$ as well as a $\phi$-measure on $\Xhat$ and a
hermitian metric on $\Ehat$ so that
\begin{enumerate}
 \item[(a)] $\Phat$ restricts to $X$ in the sense of Definition \ref{def:extension operator}
     and $\Phat_X$ is a $\phi$-differential operator on $X$ augmenting $P$,
 \item[(b)]
 $\ker \Phat_X \cap \ker \gamma$ is a finite dimensional subspace of $\CdotsBCinf(X;\Ehat)$,
 \item[(c)] $\Phat+\Pi_{sh}$ is fully elliptic, self-adjoint and invertible in $L^2_\phi(\Xhat;\Ehat)$, where $\Pi_{sh}$ is the orthogonal projection from $L^2_\phi(\Xhat;\Ehat)$ to the space of shadow solutions in (b).
\end{enumerate}
\end{theorem}
In (b) the kernel of $\Phat_X$ is taken in $\calA_\s(X)$.
In (c) we consider $\CdotsBCinf(X;\Ehat)$ as a subspace of $L^2_\phi(\Xhat;\Ehat)$ as in \eqref{eqn:CdotBCinf identific}. For $\Xhat$ one can always take the double of $X$, but other extensions are possible.

In the case that $P=-\Delta_g$ is the Laplacian of a $\phi$-metric $g$ there is a quite straight-forward proof of Theorem \ref{thm:construction inv ext}, which we present in Subsection \ref{ssec:laplacian.extension}. It motivates the last steps (achieving full ellipticity, Subsection \ref{ssec:ext.full}) of the general construction, which  is much more involved.

We now give an outline of the general construction of $\Phat$. It proceeds in several steps, summarized as follows:
\newcommand{\rsa}{\rightsquigarrow}
\begin{equation}\label{eqn:inv ext outline}
 \newcommand{\move}[1]{\hspace{-#1pt}} \newcommand{\nxt}[1][3]{\move{#1}&\move{#1}}
\begin{array}{ccccccccccccccc}
 P \nxt[5] \stackrel{\text{aug}}{\rsa} & \Pbar & \stackrel{\text{ext}}\rsa \nxt \Phat_0
  \nxt \stackrel{\text{aug}}\rsa \nxt \Phatbar_0 \nxt \stackrel{\text{add}}\rsa
  \nxt[5] \Phat_1 \nxt[5] \stackrel{\text{aug}}\rsa \nxt \Phatbar_1
  \nxt \stackrel{\text{add}}\rsa \nxt[5] \Phat \nxt \stackrel{\text{mod}}\rsa \nxt[1] \Phat + \Pi_{sh} \\[2pt]
 \nxt[5] \nxt\text{\small s.a.} \nxt \nxt \text{\small $\phi$-ell} \nxt \nxt[5] \text{\small + s.a.}
  \nxt[5] \nxt[5] \text{\small fully ell} \nxt[5] \nxt \text{\small + s.a.} \move{5} & \nxt \text{\small + supp${}^+$} \nxt \nxt \text{\small + inv} \\[2pt]
 E,E' \nxt[5] \nxt \Ebar \nxt \nxt \Ehat_0 \nxt \nxt[5] \Ehat_1 \nxt[5] \nxt[5] \Ehat_1 \nxt[5] \nxt \Ehat
  \nxt[5] \nxt \Ehat \nxt \nxt \Ehat
\end{array}
\end{equation}
Here `aug' means augmentation (enlarging the bundle), `ext' means extension (from $X$ to $\Xhat$), `add' means that a smoothing perturbation term is added that is supported in $X^-$ (i.e.\ restricts to $X$ as the zero operator) while `mod' means modification. The term `s.a.' means formally self-adjoint, `ell' means elliptic,
`supp${}^+$' means that all functions in the kernel are supported in $X$
 and `inv' means invertible.

There are three main steps here: constructing a $\phi$-elliptic extension $\Phat_0$, achieving full ellipticity with $\Phat_1$ and then achieving invertibility of $\Phat+\Pi_{sh}$.
In each of these constructions it is useful to start with a formally self-adjoint operator, so each step is preceded by an augmentation whose only purpose is to make the previous operator formally self-adjoint; this is achieved by considering the augmentation $\Pbar=
\big(\begin{smallmatrix}
 0 & P^\star \\ P & 0
\end{smallmatrix}\big)
$. For this purpose auxiliary metrics will be chosen.
These augmentations and the vector bundles (third line in \eqref{eqn:inv ext outline}) are introduced in Subsection \ref{ssec:inv ext bundles aug}.

The construction of $\Phat_0$ from $P$ via $\Pbar$
 is  the same as done by Seeley, adapted to the $\phi$-setting, cf.\ (i') in Section \ref{subsec:inv ext}, and also the step from full ellipticity to invertibility ($\Phat_1$ to $\Phat+\Pi_{sh}$) is analogous to Seeley's construction, except for a slight simplification afforded by constructing $\Phatbar_1$ first. We do this in Subsections \ref{ssec:extension phi-elliptic} and \ref{ssec:ext.inv}.

 The main new contribution in this paper is the construction of $\Phat_1$, carried out in Subsection \ref{ssec:ext.full}. We now give an outline of this construction.
 Recall that the normal family of $\Phatbar_0$ consists of $\Psi$DOs $N(\Phatbar_0)(\mu)$ on the fibres $F_y$, where $\mu=(\tau;y,\eta)\in\R\times T^*B$. So the task is to perturb this normal family to make it invertible.
The first idea would be to proceed as in step (iii') in Section \ref{subsec:inv ext}, i.e.\
to add a projection to a subspace $V_\mu$ which is complementary to the range of $N(\Phatbar_0)(\mu)$, for each $\mu$. However, the dimension of the kernels of $N(\Phatbar_0)(\mu)$
may vary with $\mu$, so there is no continuous family of such projections. We circumvent this problem by adding $i=\sqrt{-1}$ times an orthogonal projection, which yields invertibility under the weaker condition that $V_\mu+\rg N(\Phatbar_0)(\mu)$ is the full space (the sum not necessarily being direct), see Lemma \ref{lem:proj.inversion}. Construction of a smooth finite rank bundle $(V_\mu)_{\mu\in\R\times T^*B}$ satisfying this condition is possible since $N(\Phatbar_0)(\mu)$ is invertible for $\mu$ outside a compact set.
Then we obtain $\Phat_1$ from $\Phatbar_0$ by adding a smoothing operator
whose normal family at $\mu$ is 
$i$ times the orthogonal projection to $V_\mu$.
In order to ensure that $\Phat_1$ still extends the original operator $P$, we need to construct $V_\mu$ so that its elements are supported in $X^-$. This is where Assumption \ref{UCNF} is needed, see Lemma \ref{lem:complement in M-}. See also Remark \ref{rem:need UCNF}.

We remark that (a) in Theorem \ref{thm:construction inv ext} could be strengthened to $\Phat$ being a $\phi$-differential operator \textit{near} $X$, in the sense that
its Schwartz kernel has support in $\diag_\Xhat \cup \interior{X^-}^2$.
This follows from the fact that the space $W$ constructed in Lemma \ref{lem:complement in M-}(ii) is actually contained in $\CdotsBCinf(X^-;\Ehat)$. However, we don't need this extra information, and our statement seems cleaner.

\smallskip

Throughout this section we fix a background $\phi$-metric $\ghat$ on $\Xhat$, with respect to the extended fibration $\phihat$ of $\partial \Xhat$, and also hermitian metrics on the vector bundles $E,E'$.

\subsection{The unique continuation assumption}\label{ssec:ucp asssumption}
Let us say that a differential operator $T$ on a compact manifold $F$ with (\textBC-)boundary has the \textbf{unique continuation property (UCP) with respect to the boundary} if
\begin{equation}
 \label{eqn:UCP def}
 \ker T \cap \Cdotinf(F,E) = \{0\} \,.
\end{equation}
That is, any solution $u$ of $Tu=0$ vanishing to infinite order at $\partial F$ must vanish identically. If $T$ is elliptic of order $m$ then this is equivalent to $Tu=0, \gamma u =0\Rightarrow u=0$. We also express this by saying that \textbf{$T$ has no shadow solutions.}

Recall  that
the normal family of $P$ is a family of differential operators on $F$,
$N(P)(\mu):\Cinf(F_y,E_y) \to \Cinf(F_y,E_y')$, and that
 $N(P^\star)(\mu)=N(P)(\mu)^\star$ (see Appendix \ref{sec:app phi ops}).

\begin{assumption}\label{UCNF} We assume that the normal families of
 $P$ and $P^\star$ have the unique continuation property with respect to the boundary of $F$.
\end{assumption}

See Remark \ref{rem:need UCNF} for some considerations on this assumption.
Note that to define the adjoint we need to choose a $\phi$-metric on $X$ and bundle metrics on $E,E'$. However, the assumption on $N(P)^\star$ is  independent of the choice of metrics on $X$ and $E, E'$. This can be proved as follows. First, by considering $T=N(P)(\mu)$ and $T^\star$ as operators $V \to V'$ resp. $V'\to V$ where $V=\CdotBCinf(F_y,E_y)$, $V'=\CdotBCinf(F_y,E_y')$ and using $\ker T^\star = (\rg T)^\perp$
we see that the failure of the UCP for $T^\star$ is equivalent to the existence of $v\in V'\setminus\{0\}$ such that $v\perp \rg T$. Here we still use the $L^2$ scalar product in $V'$. Next, changing the metrics on $X$ and $E'$ amounts to replacing the $L^2$-scalar product $\langle\ ,\ \rangle$ on $V'$ by $\langle u,v\rangle_1 = \langle u,Av\rangle$ where $A$ is a bundle automorphism of $E'$. This implies that $v\perp_1 \rg T \iff Av\perp \rg T$. Now $v\mapsto Av$ is an invertible linear map $V'\to V'$, and the claim follows.

\subsection{The case of the Laplacian}\label{ssec:laplacian.extension}
If $P=-\Delta_g$ is the scalar Laplacian of a $\phi$-metric $g$ then there is a much simpler proof of Theorem \ref{thm:construction inv ext} than in the general case.  We present this first as it motivates the last steps  
of the general construction, i.e. the construction of $\Phat_1$ and of $\Phat$. The other steps are trivial in this case. The bundles are all trivial line bundles, but the same construction also works for the Hodge Laplacian on forms.

First, $P$ has a natural extension to $\Xhat$ as $-\Delta_{\ghat}$ for a $\phi$-metric $\ghat$ extending $g$. This is a $\phi$-elliptic, self-adjoint differential operator on $\Xhat$. However, it is not fully elliptic since its normal family is $N(-\Delta_{\ghat}) (\tau;y,\eta) = - \Delta_{\Fhat_y} + |(\tau,\eta)|^2$, with the induced metric $|\cdot|$ on $(\tau,\eta)$ space, and for  $(\tau,\eta)=0$ this operator vanishes on the constants.
This can be remedied easily with the help of the positivity of the Laplacian: Choose a function $a \in \Cinf(\Xhat, \R)$ satisfying
\begin{enumerate}
 \item\label{lapl.ext.item1} $a=0$ on $X$,
 \item\label{lapl.ext.item2} $a\geq 0$ everywhere on $\Xhat$ and
 \item\label{lapl.ext.item3} $a_{|\Fhat_y}$ does not vanish identically for any $y \in B$,
\end{enumerate}
and set $\Phat = -\Delta_{\ghat} + a$. This is a $\phi$-differential operator whose normal family is
$$ N(\Phat) (\tau;y,\eta) = - \Delta_{\Fhat_y} + |(\tau,\eta)|^2 + a_{|\Fhat_y} \in \Diff^2(\Fhat_y)$$
which is non-negative by \eqref{lapl.ext.item2} and even positive by \eqref{lapl.ext.item3}, for each $ (\tau;y,\eta)\in \R \times T^*B$, so $\Phat$ is fully elliptic. Also, $\Phat$ is positive for the same reason, so it is invertible. Finally, $\Phat$ extends $P$ by \eqref{lapl.ext.item1}, so it satisfies the claims of Theorem \ref{thm:construction inv ext}, where $\ker \Phat_X=\{0\}$ and hence $\Pi_{sh}=0$ since the Laplacian has the unique continuation property.

\subsection{The augmentations and vector bundles}\label{ssec:inv ext bundles aug}
We first discuss the three augmentation steps in \eqref{eqn:inv ext outline} and define the vector bundles on which the operators act. Each operator in \eqref{eqn:inv ext outline} except $P$ acts from sections of the bundle noted underneath it to sections of the same bundle. Recall that we consider $P$ as an operator from $\calF(X;E)$ to $\calF(X;E^\prime)$.

First, we let
 \begin{equation}\label{eqn:def.Pbar}
  \Pbar =
   \begin{pmatrix}
    0 & P^\star\\ P & 0
   \end{pmatrix} : \calF(X;\Ebar)\to\calF(X;\Ebar),\quad \Ebar=E\oplus E'
 \end{equation}
where $P^\star$ is the formal adjoint of $P$ with respect to the chosen $\phi$-metric on $X$ and hermitian metrics on $E,E'$. Thus, $\Pbar$ is an augmentation of $P$ with respect to the bundle maps
\begin{equation}
 \label{eqn:def i_0 pi_0}
  E \stackrel{\iota_0}\longhookrightarrow \Ebar \stackrel{\pi_0}\longtwoheadrightarrow E,
\quad
  E' \stackrel{\iota_0'}\longhookrightarrow \Ebar
 \stackrel{\pi_0'}\longtwoheadrightarrow E'
\end{equation}
which are injection as and projection to the first and second factor.
Next, we extend $\Ebar\to X$ to $\Ehat_0\to\Xhat$ using the extensions of $E,E'$ to $\Xhat$, which exist by the assumptions of Theorem \ref{thm:construction inv ext}.
For $j=0,1$ and assuming that $\Phat_j$ is already constructed,%
\footnote{$\Phat_0$ and $\Phat_1$ will be constructed in Section \ref{ssec:extension phi-elliptic}.} we now let
 \begin{equation}\label{eqn:def.Phatbar}
  \Phatbar_j =
   \begin{pmatrix}
    0 & \Phat_j^\star\\ \Phat_j & 0
   \end{pmatrix} : \calF(\Xhat;\Ehat_{j+1})\to\calF(\Xhat;\Ehat_{j+1}),\quad \Ehat_{j+1}=\Ehat_j\oplus \Ehat_j
 \end{equation}
where
\begin{equation}
 \label{eqn:def i_j pi_j}
  \Ehat_j \stackrel{\iota_{j+1}}\longhookrightarrow \Ehat_{j+1} \stackrel{\pi_{j+1}}\longtwoheadrightarrow \Ehat_j,
\quad
  \Ehat_j' \stackrel{\iota_{j+1}'}\longhookrightarrow \Ehat_{j+1} \stackrel{\pi_{j+1}'}\longtwoheadrightarrow \Ehat_j'
\end{equation}
are injection as and projection to the first and second factor. We also write $\Ehat=\Ehat_2$.

Recall that our notion of restriction, Definition \ref{def:extension operator}, is preserved under taking adjoints. Therefore, if $\Phat_j$ restricts to $X$ then so does $\Phatbar_j$, and if $(\Phat_j)_X$ is a $\phi$-differential operator then so is $(\Phatbar_j)_X$.

On $\Ebar$  we use the hermitian product metric, starting with the given metrics on $E$ and $E'$. The metric on $\Ehat_0$ is any extension of this metric. On $\Ehat_1$, $\Ehat=\Ehat_2$ we use the hermitian product metrics.
Then $\Pbar$, $\Phatbar_0$, $\Phatbar_1$ are formally self-adjoint $\phi$-elliptic operators.

Over $X$, the bundles $E,E'$ are related to $\Ehat$ via the compositions
\begin{equation}
 \label{eqn:relation E Ehat}
  E \stackrel{\iota}\longhookrightarrow \Ehat_X \stackrel{\pi}\longtwoheadrightarrow E\,,\qquad \iota = \iota_2\iota_1\iota_0\,,\quad \pi = \pi_0\pi_1\pi_2
\end{equation}
and similarly $E' \stackrel{\iota'}\longhookrightarrow \Ehat_X
 \stackrel{\pi'}\longtwoheadrightarrow E'$ where
 $\iota' = \iota_2'\iota_1'\iota_0'$, $\pi' = \pi_0'\pi_1'\pi_2'$.
The restrictions $(\Phatbar_1)_X$, $(\Phatbar_2)_X$ could be represented as $4\times 4$ and $8\times8$ matrices whose only non-zero entries are an alternating sequence of $P$ and $P^\star$ on the antidiagonal. The 'original' $P$ is included as lower left corner, or formally:
\begin{equation}
 \label{eqn:inv ext P in P2}
  P = \pi' \,(\Phatbar_2)_X \,\iota\,.
\end{equation}

\subsection{Two functional analytic lemmata}
The following elementary facts are used at several places in the construction.

\begin{lemma}\label{lem:proj.inversion} Let $H$ be a complex Hilbert space, $T : \mathrm{dom}\, T \subset H \to H$ a densely defined self-adjoint Fredholm operator and $\Pi$ a finite rank orthogonal projection in $H$.
Let $\alpha>0$. Then:
\begin{enumerate}
 \item[(a)]
  $T+\alpha\Pi$ is invertible if $\rg T\oplus\rg \Pi=H$.
 \item[(b)]
  $T + i\alpha\Pi$ is invertible if and only if $\rg T + \rg \Pi = H$.
\end{enumerate}
\end{lemma}
Statement (a) has a weak converse: $T+\alpha\Pi$ invertible implies $\rg T + \rg \Pi = H$ (same proof as for (b)), but the sum need not be direct as the example $T=\Pi=\id$ on a finite dimensional space shows.
\begin{proof}
The operators in (a) and (b)
are Fredholm with index zero since $T$ has this property
and $\Pi$ has finite rank. So for invertibility it suffices to check injectivity in each case.
Also, observe that $\rg T+\rg \Pi=H$ implies $\ker T\cap\ker \Pi=\{0\}$ since
$(\ker T \cap \ker \Pi)^\perp = (\ker T)^\perp + (\ker \Pi)^\perp = \rg T + \rg\Pi$.

(a) If $(T+\alpha\Pi)u=0$ then $Tu=-\alpha\Pi u$, so $Tu=\Pi u=0$ since
$\rg T\cap\rg \Pi=\{0\}$. Then $\ker T\cap\ker \Pi=\{0\}$ implies $u=0$.

(b) If $(T + i\alpha\Pi)u=0$ then $\langle Tu,u\rangle + i\alpha\langle \Pi u,u\rangle=0$, and because both $T$ and $\Pi$ are self-adjoint this implies $\langle \Pi u, u\rangle = 0$, hence $\Pi u = 0$ since $\Pi$ is an orthogonal projection. Inserting this into $(T+\alpha i\Pi)u=0$ we get $T u = 0$, and
again $\ker T\cap\ker \Pi=\{0\}$ implies $u=0$.

For the converse suppose that $T + i\alpha\Pi$ is invertible. Given $g \in H$, there is $f \in H$ so that $(T + i\alpha\Pi) f = g$. Then $g = Tf+\Pi(i\alpha f)\in \rg T+\rg\Pi$.
\end{proof}

The following lemma allows us to find projections as in Lemma \ref{lem:proj.inversion}.

\begin{lemma} \label{lem:complement in M-}
 Let $M$ be a manifold with \textBC-boundary and $\Mhat$ be a \textBC-extension of $M$.

 Then for any finite-dimensional subspace $K\subset \Cdotsinf(\Mhat)$ the following are equivalent:
  \begin{enumerate}
   \item[(i)] $ u\in K,\ \supp u \subset M \Rightarrow u = 0$
   \item[(ii)] There is a subspace $W \subset \CdotsBCinf(M^-)$ satisfying
       $$ W \oplus K^\perp = L^2(\Mhat)\,. $$
  \end{enumerate}
 The analogous statements hold for sections of a hermitian vector bundle over $\Mhat$.
\end{lemma}
Here we use the $L^2$ scalar product defined by any conormal density on $\Mhat$.
The lemma says that if all non-trivial elements of $K$ are already non-trivial on $M^-$ then the complement $K$ of $K^\perp$ may be replaced by a space of sections supported in $M^-$. Note that condition (i) may be expressed as '$u\in K$ is determined by $u_{|M^-}$' and therefore is a unique continuation condition.

In the applications of the lemma $K$ is the kernel of a self-adjoint ($\phi$-)elliptic $\Psi$DO $\That$ on $\Mhat$, and then $W$ is a complement to its range.
We need this lemma in two settings: In the first setting $\That$ is an operator in the normal family of $\Phatbar_0$ and $M$ is a fibre $F$ (proof of Lemma \ref{lem:bundle.ext}). Here $\dsM=\emptyset$. In the second setting $\That$ is a modification of $\Phatbar_1$ and $M=X$ (see the proof of Theorem \ref{thm:construction inv ext} in Subsection \ref{ssec:ext.inv}). We will not need the analogous statement for general function spaces $\calF$ though.

The lemma is a smooth version of a simple fact about Hilbert spaces:
Let $H=L^2(\Mhat)$, $H^\pm=L^2(M^\pm)$. Then $H = H^+ \oplus H^-$ (orthogonal direct sum), and for a closed subspace $K \subset H$ we have from $H^-=(H^+)^\perp$ that
$$ K\cap H^+ = \{0\} \iff K^\perp + H^- = H \iff \exists W\subset H^-:\ K^\perp \oplus W = H\,.$$

\begin{proof}
(ii)$\Rightarrow$(i): If $W\subset \CdotsBCinf(M^-)$ satisfies (ii) and $u\in K$ is supported in $M$ then $u\perp W$ and $u\perp K^\perp$, so $u\perp W\oplus K^\perp=L^2(\Mhat)$, hence $u=0$.

(i)$\Rightarrow$(ii):
Let $\rho:\Mhat\to\R$ be a defining function for $\partial_{\BC}M$ so that $M=\{\rho\geq0\}$. We first show that we can enlarge $M$ slightly in condition (i), that is, there is $\eps>0$ so that
\begin{equation}
 \label{eqn:cond(i)enlarged}
u\in K, \supp u \subset \{\rho\geq-\eps\} \Rightarrow u=0\,.
\end{equation}
To show this, assume it was wrong. Then we could find a sequence $u_m\in K$ with $\supp u_m\subset\{\rho\geq-\frac1m\}$ and $\|u_m\|=1$. Since $\dim K<\infty$, there would be a convergent subsequence $u_{m'}\to u$ with  $u\in K$, $\|u\|=1$ by compactness, with convergence in $C(\Mhat)$, so $\supp u\subset\{\rho\geq0\}=M$. This would contradict assumption (i).

Now choose $\chi\in\Cinf(\Mhat,\R)$ supported in $M^-$ and equal to 1 on $\{\rho\leq-\eps\}$. Then \eqref{eqn:cond(i)enlarged} implies that the maps $u\mapsto \chi u$, $u\mapsto \chi^2 u$ are injective on $K$. Let $W=\chi^2 K$. We claim that $W$ satisfies $W \oplus K^\perp = L^2(\Mhat)$. To prove this, it suffices to check $W \cap K^\perp = \{0\}$ because $\dim W=\dim K$.
Now if $w\in W$ then $w=\chi^2 u$ with $u\in K$, so if $w\in K^\perp$ also then $0=\langle w,u\rangle=\langle \chi^2 u,u\rangle = \langle \chi u, \chi u\rangle$, so $\chi u=0$, hence $u=0$.
\end{proof}

\subsection{Constructing a \texorpdfstring{$\phi$}{\straightphi}-elliptic augmented extension}\label{ssec:extension phi-elliptic}
Here we construct the operator $\Phat_0$ in \eqref{eqn:inv ext outline}. Recall from \eqref{eqn:def.Pbar} that we already constructed $\Pbar$.

\begin{proposition}\label{prop:ext1}
Let $\Pbar\in\Diff^m_\phi(X;\Ebar)$ be formally self-adjoint and $\phi$-elliptic. Let $\Xhat$ be a $\phi$-\textBC-extension of $X$ and $\Ehat_0\to\Xhat$ an extension of $\Ebar$.
Then there is a $\phi$-elliptic extension $\Phat_0\in\Psi^m_\phi(\Xhat;\Ehat_0)$ of $\Pbar$.
\end{proposition}

\begin{proof}
Extend the hermitian metric on $\Ebar$ to $\Ehat_0$.
Since $\Pbar$ is assumed to have coefficients smooth up to $\dBCX$, we may extend it to a neighborhood $\tilde X$ of $X$ in $\Xhat$ as a formally self-adjoint $\phi$-elliptic differential operator $\Ptilde$.

Since the $\phi$-principal symbol $\SymP = \phisigma_m(\Ptilde)$ is $\phi$-elliptic and 
self-adjoint over $\tilde{X}$, we can extend it to a $\phi$-elliptic symbol $\SymP_1$ of order $m$ over $\Xhat$ by setting $\SymP_1=|\xi|^m\Id_{\Ehat_0}$ outside a neighborhood of the closure of $\tilde{X}$ and connecting endomorphisms of fibres of $\Ehat_0$ to the identity by moving their spectrum through the upper complex half plane, cf.\ \cite[p.\ 298]{See69}.

We choose an operator
 \[\text{$\Ptilde' \in \Psi_\phi^m(\Xhat;\Ehat_0)$ so that $\phisigma_m(\Ptilde') = \SymP_1$.} \]
This operator need not extend $\Pbar$ although its principal symbol extends the principal symbol of $\Pbar$. Therefore, we 'glue' $\Ptilde'$ away from $X$ with $\Ptilde$ on $\Xtilde$, by choosing cut-off functions
$\psi_+$, $\psi_- \in \Cinf(\Xhat,\R_+)$ so that $\psi_+ + \psi_- = 1$ on $\Xhat$, $\psi_- = 0$ in a neighborhood of $X$ and $\psi_- = 1$ outside of $\tilde{X}$, and setting
 \[ \Phat_0 = \psi_+ \Ptilde \psi_+ + \psi_- \Ptilde' \psi_- \,. \]
Then $\Phat_0 \in \Psi^m_\phi(\Xhat;\Ehat_0)$ is $\phi$-elliptic and an extension of $\Pbar$ in the sense of Definition \ref{def:extension operator}.
\end{proof}
The extension $\Phat_0$ cannot in general be chosen to be a differential operator.
Also, the construction yields a non-selfadjoint operator in general since the extended symbol $\SymP_1$ cannot in general be chosen to be self-adjoint; the reason for this is that the set of self-adjoint invertible $N\times N$ matrices is not connected for any $N$.

\subsection{Achieving full ellipticity}\label{ssec:ext.full}
We now construct the fully elliptic operator $\Phat_1$ from $\Phatbar_0$ in \eqref{eqn:inv ext outline}, where $\Phatbar_0$ is an augmentation  as in \eqref{eqn:def.Phatbar} of the operator $\Phat_0$ just constructed.
For this we need to perturb the normal family, and it is in this step that we need to make the Assumption \ref{UCNF} on our original operator $P$.

Assumption \ref{UCNF} implies that the normal family of $(\Phatbar_0)_X$ also has the unique continuation property with respect to the boundary since this operator is a direct sum of two copies of $P$ and $P^\star$ each. Therefore, we need to prove the following proposition.

\begin{proposition}\label{prop:construction fully elliptic}
Let $\Phatbar_0\in\Psi_\phi^m(\Xhat;\Ehat_1)$ be $\phi$-elliptic and selfadjoint. Also, assume that it restricts to $X$ (Definition \ref{def:extension operator}) and that the normal family of  $(\Phatbar_0)_X$ has the unique continuation property with respect to the boundary (as defined before Assumption \ref{UCNF}).

Then there is a \emph{fully} elliptic $\Phat_1\in\Psi_\phi^m(\Xhat;\Ehat_1)$ which restricts to $X$ and so that
\begin{gather}
 \label{eqn:Phat1 and Phatbar0}
 (\Phat_1)_X = (\Phatbar_0)_X\,,\quad \Phat_1 - \Phatbar_0 \in \Psi^{-\infty}_\phi (\Xhat;\Ehat_1)
\end{gather}
\end{proposition}
Note that Definition \ref{def:extension operator} involves an implicit choice of function space $\calF$ but following Remark \ref{rmk:extension operator}, this choice does not matter for $\phi$-$\Psi$DOs, and that \eqref{eqn:Phat1 and Phatbar0} is equivalent to the Schwartz kernel of $\Phat_1-\Phatbar_0$ being a smooth section on $\Xhat^2_\phi$ supported in $(X^-)^2_\phi$.

For the proof of the proposition we need two lemmata.
\begin{lemma}\label{lem:ind.sect} Let $\calV \to Z$ be a smooth vector bundle of infinite rank over a manifold $Z$ and let $w_1, \dotsc, w_k$ be smooth sections of $\calV$. Then there are smooth sections $s_1, \dotsc, s_k$ of $\calV$ so that, for all $\eps \neq 0$,
 \[ w_1 + \eps s_1, \dotsc, w_k + \eps s_k \]
are linearly independent sections of $\calV$.
\end{lemma}

\begin{proof} For $z \in Z$, let $W(z) = \{ w_1(z), \dotsc, w_k(z)\}$. As both the base $Z$ and $\spn W(z)$ are finite dimensional, we can choose a section $s_1$ of $\calV$ so that $s_1(z) \not\in \spn W(z)$ for all $z \in Z$. Step by step, choose sections $s_1, \dotsc, s_k$ of $\calV$ so that
 \begin{equation}\label{eqn:ind.sect.1}
  s_j(z) \not\in \spn \big( W(z) \cup S_{j-1}(z) \big) \quad\text{for all $z \in Z$,}
 \end{equation}
 where $S_l(z) = \{s_1(z), \dotsc, s_l(z) \}$.
 Fix $z$. Then inductively  $\spn W(z) \cap \spn S_l(z)=\{0\}$ for each $l$.
Then if $\eps \neq 0$ and $\sum_i \lambda_i \big(w_i(z) + \eps s_i(z)\big) = 0$
it follows that $\sum_i \lambda_i  s_i(z) = -\frac1\eps \sum_i \lambda_i w_i(z)$, so both sides must be zero, hence $\lambda_i=0$ for all $i$ since the $s_i(z)$ are linearly independent.
So we even get that $w_1(z) + \eps s_1(z), \dotsc, w_k(z) + \eps s_k(z)$ are linearly independent for each $z$.
\end{proof}

\begin{lemma}\label{lem:bundle.ext}
In the setting of Proposition \ref{prop:construction fully elliptic}  there is a smooth finite rank vector bundle $V \to \R \times T^*B$  so that for each $\mu= (\tau;y,\eta) \in \R \times T^*B$:
  \begin{gather}
  \label{eqn:V in F-}
  V_\mu \subset \CdotBCinf(F^-_y,\Ehat_1) \\
  \label{eqn:bundle.ext.2}
   V_{\mu} + \rg N(\Phatbar_0) (\mu)
    = L^2(\Fhat_y,\Ehat_{1y})\,.
  \end{gather}
\end{lemma}

\begin{proof}
In this proof write $T=\Phatbar_0$ and leave out bundles from the notation. Since $T$ is $\phi$-elliptic there is a compact set
$\calK \subset \R \times T^*B$ so that $N(T)(\mu)$ is invertible for $\mu\not\in\calK$,
see \cite[Prop.\ 2]{Mel95a}.

Fix $\mu=(\tau;y,\eta)\in\calK$. We apply Lemma \ref{lem:complement in M-} with $\Mhat=\Fhat_y$ and $K=\ker N(T)(\mu)$. Condition (i) in the lemma is simply the unique continuation property with respect to the \textBC-boundary, so it is satisfied by assumption. Now $N(T)(\mu)$ is an elliptic $\Psi$DO on the closed manifold $\Fhat_y$, hence Fredholm in $L^2(\Fhat_y)$, so $\dim K<\infty$ and $K^\perp = \rg N(T)(\mu)$ by self-adjointness, so (ii) of the lemma gives a subspace
$$ V_\mu' \subset\CdotBCinf(F^-_y) \text{ so that } V_\mu' \oplus \rg N(T)(\mu) = L^2(\Fhat_y)\,. $$
Then by continuity there is an open neighborhood $U_\mu$ of $\mu$ in $\R\times T^*B$ so that
\begin{equation}
\label{eqn:V'compl local}
V_\mu' + \rg N(T)(\tilde\mu) = L^2(\Fhat_\ytilde)
\end{equation}
for all $\tilde\mu=(\tautilde;\ytilde,\etatilde)\in U_\mu$,
 where we identify nearby fibres $\Fhat_\ytilde$.

We now combine the $V_\mu'$ using a compactness argument and Lemma \ref{lem:ind.sect} to obtain $V$:
As $\calK$ is compact, finitely many of the $U_\mu$ suffice to cover this set, say $U_1, \dotsc, U_N$, where $U_j = U_{\mu_j}$. Let
$\psi_1, \dotsc, \psi_N\in\Cinf(\R\times T^*B)$ be a partition of unity subordinate to the cover $U_1, \dotsc, U_N$ of $\calK$.
Choose a basis $(\omega_{ij}')_i$ for each $V_{\mu_j}'$ and apply Lemma \ref{lem:ind.sect} to the set of $\omega_{ij} = \psi_j \omega_{ij}'$ over all $i,j$,
where $Z=\R\times T^*B$ and the bundle $\calV\to\R\times T^*B$ is given by $\calV_\mu = \CdotBCinf(F^-_y)$ for $\mu=(\tau;y,\eta)$. The lemma gives sections $s_{ij}$ of $\calV$ so that for each $\eps>0$ the sections $(\omega_{ij}+\eps s_{ij})_{i,j}$ are linearly independent, so span a subbundle $V^\eps$ of $\calV$.
Moreover, for each $j$ condition \eqref{eqn:V'compl local}, with $\mu=\mu_j$ and $\mutilde$ in the compact set $\supp\psi_j$, is stable under small perturbations of and under enlarging $V_\mu'$, hence is satisfied with $V_\mu'$ replaced by $V^\eps$ for $\eps$ sufficiently small. Therefore we may choose $\eps$ sufficiently small so that \eqref{eqn:V'compl local}, with $V_\mu'$ replaced by $V^\eps_\mu$ and $\mutilde=\mu$, holds for all $\mu \in \calK$. Since $\rg N(T)(\mu) = L^2(\Fhat_y)$ for $\mu \not\in \calK$ the bundle $V=V^\eps$ satisfies the claim of the proposition.
\end{proof}

\begin{proof}[Proof of Proposition \ref{prop:construction fully elliptic}]
Construct a bundle $V$ as in Lemma \ref{lem:bundle.ext}. For each $\mu=(\tau;y,\eta)\in \R\times T^*B$ let $\widetilde{\Pi}_\mu$ be the orthogonal projection
 \[ \widetilde{\Pi}_\mu : L^2(\Fhat_y,\Ehat_{1y}) \longrightarrow V_{\mu} \,. \]

We now assemble the Schwartz kernels of these projections to construct the desired perturbation of $\Phatbar_0$.

Choose a compact set $\calK\subset\R\times T^*B$ so that $N(\Phatbar_0)(\mu)$ is invertible for $\mu\not\in\calK$ (see the beginning of the proof of Lemma \ref{lem:bundle.ext}) and a smooth,  compactly supported function $e : \R \times T^*B \to \R_+$ so that $e >0$ on $\calK$. Let $\kappatilde(\mu;z,z')$, $z,z'\in \Fhat_y$, be the Schwartz kernel of $\widetilde{\Pi}_\mu$. By \eqref{eqn:V in F-} it is supported in $(F^-)^2$.
Taking the $(\tau,\eta)\mapsto (T,Y)$ inverse Fourier transform of $e\,\kappatilde$ we obtain a smooth function $\kappa(T,y,Y;z,z')$ vanishing rapidly as $|(T,Y)|\to\infty$.
Recall that $(T,y,Y;z,z')$ are coordinates on the interior of the $\phi$-front face $\phi f$ of $\Xhat^2_\phi$, and the boundary of $\phi f$ is the bundle of spheres at infinity $|(T,Y)|=\infty$. Therefore $\kappa$ is smooth on $\phi f$ and vanishes to infinite order at its boundary faces, and is supported in $(z,z')\in (F^-)^2$ (see Appendix \ref{sec:app phi ops}). Hence $\kappa$ can be extended to a smooth section, which we still denote by $\kappa$, on all of $\Xhat_\phi^2$ which vanishes to infinite order at all boundary hypersurfaces except $\phi f$ and which is supported in $(X^-)^2_\phi$.
Then the operator $\widetilde{\Pi}$ with Schwartz kernel $\kappa$ is in $\Psi^{-\infty}_\phi(\Xhat;\Ehat_1)$ and has normal family
$$ N(\widetilde{\Pi})(\mu) = e(\mu)\widetilde{\Pi}_\mu$$
by construction.

Define $\Phat_1=\Phatbar_0+i\widetilde{\Pi}$.
As ${\widetilde\Pi}$ is of order $-\infty$, $\Phat_1$ is still $\phi$-elliptic. Moreover,
$N(\Phat_1)(\mu)=N(\Phatbar_0)(\mu) + i e(\mu)\widetilde{\Pi}_\mu$ is invertible for all $\mu$ by  Lemma \ref{lem:proj.inversion}(b) because of \eqref{eqn:bundle.ext.2} and since
$e(\mu)>0$ for all $\mu$ where  $N(\Phatbar_0)(\mu)$ is not invertible. By construction, $\widetilde{\Pi}$ is supported in $X^-$ so we obtain \eqref{eqn:Phat1 and Phatbar0}.
\end{proof}

\subsection{Achieving invertibility; proof of Theorem \ref{thm:construction inv ext}}\label{ssec:ext.inv}
We now prove Theorem \ref{thm:construction inv ext}.
We follow the steps outlined in and around equation \eqref{eqn:inv ext outline}: We  choose a $\phi$-metric $g$ on $X$ and hermitian metrics on  $E$, $E'$. Then we augment $P$ to $\Pbar=
\big(\begin{smallmatrix}
 0&P^\star\\ P & 0
\end{smallmatrix}\big)
$
on $\Ebar\to X$, where $P^\star$ is the formal adjoint and $\Ebar=E\oplus E'$ (with the product metric) as explained in Subsection \ref{ssec:inv ext bundles aug}. We extend $X$ to $\Xhat$, $\Ebar$ to $\Ehat_0$, and also the metrics. Then we use Proposition \ref{prop:ext1} to extend $\Pbar$ to a $\phi$-elliptic operator $\Phat_0\in\Psi^m_\phi(\Xhat;\Ehat_0)$. We augment $\Phat_0$, $\Ehat_0$ (as explained in Subsection  \ref{ssec:inv ext bundles aug}) to obtain the self-adjoint, $\phi$-elliptic operator $\Phatbar_0\in\Psi^m_\phi(\Xhat;\Ehat_1)$. Then $\Phatbar_0$ restricts to $X$ as the 
augmentation of the $\phi$-differential operator $\Pbar$ which in turn augments $P$. Then we use Proposition \ref{prop:construction fully elliptic} to find a fully elliptic $\Phat_1\in\Psi^m_\phi(\Xhat;\Ehat_1)$ which agrees with $\Phatbar_0$ on $X$.
In this step we use the unique continuation property with respect to the boundary for $N(P)$ and $N(P)^\star$, Assumption \ref{UCNF}. Now we augment $\Phat_1$, $\Ehat_1$ again to obtain $\Phatbar_1\in \Psi^m_\phi(\Xhat;\Ehat)$ where $\Ehat = \Ehat_1\oplus\Ehat_1$. Then $\Phatbar_1$ is self-adjoint, fully elliptic and restricts to $X$ as a $\phi$-differential operator which augments $P$.

It remains to make the operator invertible, i.e.\ to construct $\Phat + \Pi_{sh}$ as described in Theorem \ref{thm:construction inv ext} from $\Phatbar_1$. First we recall from Appendix \ref{sec:app phi ops} that
$$\Phatbar_1 \text{ fully elliptic} \,\Rightarrow\, \ker \Phatbar_1 \subset  \Cdotsinf(\Xhat;\Ehat)\,,$$
the space of smooth sections vanishing to infinite order at the (singular) boundary $\partial\Xhat$, and that this kernel has finite dimension.
Let $V_+=\{u\in\ker\Phatbar_1:\, \supp u \subset X\}$. Since $(\Phatbar_1)_X$ is a $\phi$-elliptic differential operator we have
\begin{equation}
 \label{eqn: ker Phatbar1}
 V_+ = \ker (\Phatbar_1)_X \cap \ker\gamma \,.
\end{equation}
Let $\Pi_{sh}$ be the orthogonal projection to $V_+$ in $L^2_\phi(\Xhat;\Ehat)$. Then $(\Phatbar_1)_X+\big(\Pi_{sh}\big)_X$ is a modification of $(\Phatbar_1)_X$ in the sense of Definition \ref{def:modification}, and the conditions around \eqref{eqn:form selfadjoint} in  Proposition \ref{prop:modification B} are satisfied, so by \eqref{eqn:modification no kernel} we have
\begin{equation}
 \label{eqn:Phatbar1 eqn1}
\ker((\Phatbar_1)_X+\big(\Pi_{sh}\big)_X) \cap \ker\gamma = \{0\}\,.
 \end{equation}
This implies that  $K := \ker(\Phatbar_1+\Pi_{sh})$ satisfies (i) of Lemma \ref{lem:complement in M-} (where $M=X$).
Also, $V_+\subset  \Cdotsinf(\Xhat;\Ehat)$ implies that $\Pi_{sh} \in x^\infty\Psi^{-\infty}_\phi(\Xhat;\Ehat)$, so $\Phatbar_1+\Pi_{sh}$ is fully elliptic, hence Fredholm, so $K$ is finite-dimensional and $ K\subset  \Cdotsinf(\Xhat;\Ehat).$
Therefore, by (ii) of Lemma \ref{lem:complement in M-} we can choose a subspace $W\subset \CdotsBCinf(X^-,\Ehat)$ complementing $K^\perp = \rg(\Phatbar_1+\Pi_{sh})$ in $L^2_\phi(\Xhat;\Ehat)$. Let $\Pi_{comp}$ be the orthogonal projection to $W$ in $L^2_\phi(\Xhat;\Ehat)$ and define
$$ \Phat = \Phatbar_1 + \Pi_{comp}\,.$$
Then $\Phat + \Pi_{sh} = (\Phatbar_1 + \Pi_{sh} ) + \Pi_{comp}$ is invertible by Lemma \ref{lem:proj.inversion} (applied with $T=\Phatbar_1+\Pi_{sh}$ and with $\Pi$ replaced by $\Pi_{comp}$), it is fully elliptic since $W\subset \CdotsBCinf(X^-,\Ehat)$ and so $\Pi_{comp}\in x^\infty\Psi^{-\infty}_\phi(\Xhat;\Ehat)$, and it is self-adjoint since $\Phatbar_1$ is self-adjoint and $\Pi_{sh}$, $\Pi_{comp}$ are orthogonal projections. Since $\Pi_{comp}$ restricts to $0$ on $X$ it is clear that the spaces of shadow solutions for $\Phat$ and $\Phatbar_1$ coincide, so (c) in Theorem \ref{thm:construction inv ext} holds.
Furthermore, (a) holds since it holds for $\Phatbar_1$ and again because $\Pi_{comp}$ restricts to zero on $X$.
Finally, $(\Phatbar_1)_X=\Phat_X$ implies that $V_+=\ker\Phat_X\cap\ker\gamma$, so (b) follows. \hfill$\square$

\begin{remark}
\label{rem:need UCNF}
We do not know if Assumption \ref{UCNF} is necessary for our theorems to hold. In our proof we need it in an essential way:
Our strategy for constructing $C$ is to make the extension $\Phat$ fully elliptic and invertible. The first obstruction to invertibility of $\Phat$ is the existence of shadow solutions for $P$. If $P$ has shadow solutions, we can deal with them by adding the projection $\Pi_{sh}$ to the space of shadow solutions and using Proposition \ref{prop:modification B}.

The second obstruction to invertibility of $\Phat$ (in the $\phi$-calculus) is that
$N(\Phat)(\mu)$ must be invertible for all $\mu$. In particular, $N(P)(\mu)$ must not have shadow solutions for all $\mu$.

Now if some $N(P)(\mu)$ had shadow solutions then we could try adding an analogous projection to it.

However, there are two problems with this:
First, while adding this projection does not alter the boundary data space of $N(P)(\mu)$ (by the same argument as for $P$), it is unclear how it affects the boundary data space of $P$ itself. Second, the space of shadow solutions is unstable under perturbations, so generically its dimension will vary with $\mu$. This means that the associated family of projections is not continuous in $\mu$, so does not define a $\Psi$DO. Note that here it is not possible to use 'too big' projections as in the proof of Proposition \ref{prop:construction fully elliptic} where we modified $N(\Phat)(\mu)$ on the minus side of $\Fhat$, since by (the proof of) Proposition \ref{prop:modification B} a projection which is too big will change the boundary data space. Also, it is not possible to use a smaller projection since then it would not remove all shadow solutions.

We also mention that an assumption analogous to Assumption \ref{UCNF} was used in
 \cite{BosFur:MIFADSFF} (see also \cite{BLZ09}) to ensure that the boundary data spaces of a family of operators vary continuously with the parameter.

Also, if $\partial F=\emptyset$, i.e.\ in the case of an interior singularity, Assumption \ref{UCNF}  is equivalent to $\ker N(P)(\mu)=\ker N(P)(\mu)^\star=0$, i.e.\ invertibility of $N(P)(\mu)$ for each $\mu$. This is a necessary condition for the invertibility on $\phi$-Sobolev spaces, and therefore very natural in the context of  the $\phi$-calculus. The same remark applies to any connected component of $\dsX$ which does not intersect $\dBCX$.

Apart from certain classes of operators (e.g. scalar second order operators), where UCP is known to hold, the UCP for any specific operator is very hard to check. It was conjectured by L. Schwartz \cite[p. 97]{Sch:EDPE} that validity of UCP for an operator on a closed manifold implies validity of UCP for its adjoint. However, this remains unsettled to this day. 
We are grateful to B. Boo{\ss}-Bavnbek for pointing out this reference to us.
\end{remark}

\subsection{Proof of Proposition \ref{prop:shadow solutions}} \label{ssec:pf prop shadow}
By Theorem \ref{thm:construction inv ext}(b) the space of shadow solutions of $\Phat_X$ is contained in $\CdotsBCinf(X;\Ehat)$. This implies the corresponding statement for $P$ since $\Phat_X$ augments $P$, so $\ker P\subset \ker \Phat_X$ with respect to the inclusion $E\hookrightarrow\Ehat$. \hfill$\square$

\section{Calder\'on Projectors}\label{sec:Calderon}

In this section we prove Theorems \ref{thm:Calderon}, \ref{thm:calderon canonical inv extension} and \ref{thm:orthogonal} and Corollary \ref{cor:L2 closure of BC spaces}, following the outline given in Section \ref{sec:outline}. The arguments in the non-singular case carry over without essential changes because of two facts: On the one hand, the Schwartz kernels of operators in the $\phi$-calculus behave in a uniform way near the $\phi$-face of $X^2_\phi$, which corresponds to the singularities, and on the other hand, the \textBC-face $\dBCX$ and the singular face $\dsX$ are transversal by assumption. Then, for instance, the transmission condition is still the correct condition to ensure that limits of restrictions (as in Proposition \ref{prop:transm.cond} below) still define $\phi$-$\Psi$DOs.

\subsection{The transmission property}\label{ssec:transmission}
We will need to generalize the transmission property of $\Psi$DOs from the classical to the singular setting. For this purpose it is useful to state it in terms of conormal distributions, whose definition is recalled in Appendix \ref{sec:app mwc}.
\begin{definition}\label{def:transm.prop}
Let $Z$ be a manifold with corners and $Y\subset Z$ an interior p-submanifold.
We say that a conormal distribution $u \in I^t_\cl(Z,Y)$  satisfies the \textbf{strong transmission condition} at $Y$ if $\tau := t + \frac14\dim Z - \frac12\codim Y$ is an integer and in its local representations \eqref{eqn:def conormal} the complete symbol satisfies \eqref{eqn:pos homogeneous} for all $\lambda\in\R\setminus\{0\}$.
\end{definition}
That is, $\lambda$ is also allowed to be negative in \eqref{eqn:pos homogeneous}.  The complete symbol depends on the coordinate system of course, but this condition is easily seen to be independent of coordinates.
The condition $\tau\in\Z$ is imposed to allow a consistent choice of powers $\lambda^{\tau-j}$ for negative $\lambda$: taking $\lambda=-1$ and using \eqref{eqn:pos homogeneous} for $\eta''$ and $-\eta''$ implies that we should have $(-1)^{2\tau}=1$, so $\tau\in\Z$.

In the case where $Z=\Xhat^2$ for a compact manifold $\Xhat$ and $Y=\diag_\Xhat$ the diagonal, distributions in $I^t_\cl(Z,Y)$ are Schwartz kernels of classical pseudodifferential operators $P\in \Psi^t(\Xhat)$. In this case it is well-known (see \cite[18.2]{Hor85}, \cite{Bou:COPVBIPT}, \cite{GruHor:TP}) that if $t \in \Z$ (note that $\tau = t$ here) then the strong transmission condition implies the transmission property for $P$ with respect to any hypersurface $W\subset\Xhat$ that divides $\Xhat$ into manifolds with boundary $X^\pm$; this says that for any $f\in\Cinf(\Xhat)$ the functions $\left[P(\chi_{X^\pm}f)\right]_{|\interior{X^\pm}}$, where $\chi_{X^\pm}$ is the characteristic function of $X^\pm$, extend smoothly to $X^\pm$.
In fact, the transmission property for $P$ for a given $W$ is equivalent to a similar condition on the full symbol of $P$, called the transmission condition, only at the conormal bundle of $W$. Our condition does not refer to an a priori choice of $W$. It could be refined to the transmission condition as in  \cite{Hor85}, but we don't need this here.

\begin{proposition}
 \label{prop:transm.cond}
 Let $Z$ be a manifold with corners and $Y\subset Z$ an interior p-submani\-fold. Furthermore, let   $H,H'\subset Z$ be interior p-hypersurfaces such that $H,H',Y$ intersect pairwise normally transversally\footnote{That is, near any point $p$ of the intersection of two of them there are adapted coordinates in which both are coordinate subspaces, and the two tangent spaces at $p$ together span $T_pZ$.}. Also, assume the conormal bundle of $H$ is orientable, allowing the choice of a 'positive side' of $H$.

Define $\tilde Z = H\cap H'$, $\tilde Y = \tilde Z \cap Y$, and assume that
\begin{equation}
 \label{eqn:diag intersection}
\text{for each $p\in \tilde Y$}, \quad  N_p^*\tilde Z\cap N_p^*Y \text{ has dimension }1
\end{equation}
Then $\tilde Y\subset\tilde Z$ is a p-submanifold. If $u\in I^t_\cl(Z,Y)$ then the restriction $u_{|H'}$ is well-defined, and if $u$
satisfies the strong transmission condition then the limit (in the sense of distributions) of $u_{|H'}$ when approaching $H$ from the positive side is well-defined and defines an element of  $I^{t+1}_\cl(\tilde Z,\tilde Y)$.
\end{proposition}

\begin{figure}[tbh]
 \begin{center}
  \includegraphics[scale=0.5]{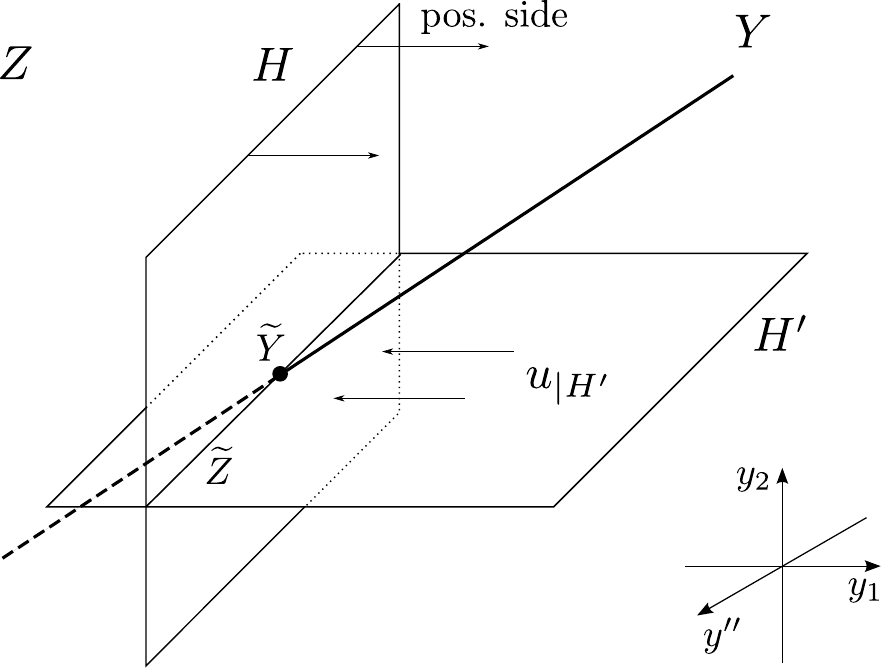}
  \caption{The p-submanifolds $Y$, $H$ and $H^\p$ of $Z$ of Prop. \ref{prop:transm.cond} in local coordinates \eqref{eqn:HH'Y coords}. Here, $\widetilde Y$ is given by the black dot that is the origin of the chosen coordinate system and the arrows next to $u_{|H^\p}$ indicate the direction $y_1 \to 0+$.}
  \label{fig:transm.cond}
 \end{center}
\end{figure}
For the last statement it clearly suffices that $u$ satisfy the strong transmission condition at $Y\cap U$ where $U$ is a neighborhood of $\tilde Z$.

In the case where $Z=\hat X^2$, $Y=\diag_{\hat X}$, $H=W\times \hat X$, $H'=\hat X\times W$ for a separating hypersurface $W\subset X$
we get $\tilde Z = W^2$, $\tilde Y = \diag_{W}$, and we recover a standard result on $\Psi$DOs closely related to the transmission property mentioned above, see  \cite[Thm.\ 18.2.17]{Hor85}, generalized below in Corollary \ref{cor:transmission property}.

\begin{proof}
Clearly the restriction and limit are well-defined outside $\tilde Y$ and the limit is smooth there, so we consider a neighborhood of a point $p\in \tilde Y$. It is easy to see
  that normal transversality and condition \eqref{eqn:diag intersection} imply that there are adapted local coordinates $x,y$ for $Z$ near $p$ in terms of which locally
\begin{equation}
\label{eqn:HH'Y coords}
 H = \{y_1=0\},\ H' = \{y_2=0\},\ Y = \{y_1=y_2, y''=0\}
\end{equation}
 where $y=(y_1,y_2,y',y'')$.
Then locally $\tilde Z= \{y_1=y_2=0\}$, $\tilde Y= \{y_1=y_2=0, y''=0\}$, in particular $\tilde Y$ is a p-submanifold of $\tilde Z$.
 The assumption $u\in I_\cl^t(Z,Y)$ means that, near $p$,
\begin{equation}
\label{eqn:TC proof}
 u(x,y_1,y_2,y',y'') = \int e^{i(y_1-y_2)\eta + iy''\eta''} a (x,y_2,y'; \eta,\eta'')\, d\eta d\eta''
\end{equation}
 where $a$ is a classical symbol of order $\tau = t + \tfrac{1}{4}\dim Z - \tfrac{1}{2}\codim Y$.

Restriction to $H'$ means setting $y_2=0$, which clearly yields a well-defined conormal distribution $u_{|y_2=0}$ on $H'$ with respect to $Y\cap H'$.
Now assume that $u$ satisfies the strong transmission condition.
We need to show that the limit of $u_{|y_2=0}$
as $y_1\to0+$ exists and defines a distribution $\tilde u$ conormal with respect to $y'' = 0$. This follows from the considerations before Theorem 18.2.17 in \cite{Hor85}: the $\eta$-integral in \eqref{eqn:TC proof} (with $y_2=0$) exists for $y_1>0$ in the sense of Lemma 18.2.16 (loc.cit.)\ as $b(x,y_1,y';\eta''):=\int^+ e^{iy_1\eta} a(x,0,y';\eta,\eta'')\,d\eta$, this is a classical symbol of order $\tau+1$ uniformly in $y_1\geq0$, and the limit $\tilde u$ of $u_{|y_2=0}$ as $y_1\to0+$ is the conormal distribution with symbol $b(x,0,y';\eta'')$. The order of $\tilde u$ is $t + 1$ since the symbol $b$ has order $\tau+1$ and $\dim \tilde Z = \dim Z - 2$, $\dim \tilde Y = \dim Y-1$, so $\frac14\dim \tilde Z - \frac12\codim \tilde Y=\frac14\dim Z - \frac12\codim Y$.
\end{proof}

\begin{corollary}
 \label{cor:transmission property}
 Let $X$ be a $\phi$-\textBC-manifold with non-empty \textBC-boundary and $\Xhat$ be an extension of $X$ across the \textBC-boundary.
 Let $Q\in\Psi^t_\phi(\Xhat)$, $t\in\Z$, and assume that the Schwartz kernel $K_Q\in I^t_\cl(\Xhat^2_\phi,\diag_\phi)$ satisfies the strong transmission condition in a neighborhood of $\diag_{\phi,\Xhat}\cap (\dBCX)^2_\phi$.

 Fix a trivialization $(-1,1)_\rho\times\dBCX$ as in \eqref{eqn:trivialize bd}.
 If $v\in\Cdotsinf(\dBCX)$ then $Q(\delta(\rho)\otimes v)$, which is smooth in $X\setminus\dBCX$, has a limit at $\dBCX$ from $\rho>0$, and  the operator
 $$ v \mapsto Q(\delta(\rho)\otimes v)_{|\dBCX} $$
thus defined is in $\Psi^{t+1}_\phi(\dBCX)$.
\end{corollary}
Note that the strong transmission condition is satisfied if $Q$ is a parametrix for a differential operator.
\begin{proof}
 Apply Proposition \ref{prop:transm.cond} with $Z=\Xhat^2_\phi$, $Y=\diag_{\phi,\Xhat}$, $u=K_Q$ and $H,H'$ the lifts to $\Xhat_\phi^2$ of $\dBCX\times\Xhat$ and $\Xhat\times\dBCX$, respectively, and the positive side of $H$ defined by $X\times\Xhat$. Transversality in the interior of $Z$ is obvious, and at the boundary of $Z$ the pairwise intersections are subsets of the interior of $\phif$, where we can use coordinates $T,x,Y,y',z,z'$ as in \eqref{eqn:phi double coords}, with $z_1$ defining $\dBCX$. In these coordinates $H=\{z_1=0\}$, $H'=\{z_1'=0\}$ and $Y=\{T=0,Y=0, z=z'\}$, so the assumptions are satisfied, with $T,Y,z_2-z_2',\dots,z_f-z_f'$ the coordinates $y''$ in \eqref{eqn:HH'Y coords}. Also, $\tilde Z = (\dBCX)^2_\phi$ and $\tilde Y$ is its $\phi$-diagonal.

 We may choose the coordinates such that $z_1=\rho$. Then the Schwartz kernel of the operator $v\mapsto Q(\delta(\rho)\otimes v)$ is $u_{|H'}$. By the proposition the limit $\tilde u$ of $u_{|H'}$ at $H$ from the positive side exists, so the limit of $Q(\delta(\rho)\otimes v)$ at $\dBCX$ from $X\setminus\dBCX$ exists, and is given (as an operator applied to $v$) by the Schwartz kernel $\tilde u$. By the proposition again, we have  
$\tilde u\in I^{t+1}_\cl((\dBCX)^2_\phi,\diag_{\phi,\dBCX})$. Also, since $u=K_Q$ vanishes to infinite order at all boundary hypersurfaces of $\Xhat^2_\phi$ except $\phif$, the analogous statement holds for $\tilde u$.
\end{proof}

\subsection{Construction of a Calder{\'o}n projector}\label{ssec:Calderon}
We now prove Theorems \ref{thm:Calderon} (i) and \ref{thm:calderon canonical inv extension} and Corollary \ref{cor:L2 closure of BC spaces}.
First, note that as the boundary data spaces $\calB_{P,\calF}$ are defined in terms of the homogeneous equation $Pu=0$, which does not see the $x^{-cm}$ factor, we may assume $c=0$.

\begin{proof}[Proof of Theorem \ref{thm:Calderon} (i)] 
Let $\Xhat$ be the \textBC-double of $X$ with extended $\phi$-structure and extensions of the bundles $E,E'$ as discussed in Subsection \ref{sssec:extspaces}.
 Let $\Ehat$, $\Phat$ and $\Pi$ be as in Theorem \ref{thm:construction inv ext}. Since $\Xhat$ has no \textBC-boundary, standard $\phi$-$\Psi$DO theory (see Appendix \ref{sec:app phi ops}) implies that $(\Phat+\Pi)^{-1}\in \Psi^{-m}_\phi(\Xhat;\Ehat)$.
We apply the construction explained in Subsection \ref{ssec:outline Calderon constr} to $\Phat+\Pi$ and then show that it yields a Calder\'on projector for $\Phat$ and then for $P$.

Choose a trivialization $(-1,1)_\rho\times\dBCX$ of a neighborhood of $\dBCX$ in $\Xhat$ as in \eqref{eqn:trivialize bd} and define $\gamma$, $\gamma^\star$ by the formulas in \eqref{eqn:cauchy.data}, \eqref{eqn:delta.data}.
Define the jump operator $\Jp_\Phat$ for $\Phat$, a differential operator of order $m-1$  from boundary data to boundary data, by the equation
 \begin{equation}\label{eqn:dfn.jump.op}
  \Phat(u^0)=(\Phat u)^0 + \gamma^\star\Jp_\Phat\gamma u
 \end{equation}
where $u$ is a function on $X$ and $u^0$ its extension to $\Xhat$ by zero.
The same equation then holds with $\Phat$ replaced by $\Phat+\Pi$ everywhere and with $\Jp_{\Phat+\Pi}=\Jp_\Phat$. This can be seen as follows. \eqref{eqn:dfn.jump.op} is equivalent to $[\Phat,\chi]=\gamma^\star\Jp_\Phat\gamma$ where  $\chi:\Xhat\to\R$ is the characteristic function of $X$. Since the Schwartz kernel of $\Pi$ is smooth and supported in $X\times X$, we have $[\Pi,\chi]=0$, so we get
 $[\Phat+\Pi,\chi]=\gamma^\star\Jp_\Phat\gamma$
also, which was to be shown.

Let
\begin{equation}
 \label{eqn:dfn.Chat}
 \Chat = \gamma (\Phat+\Pi)^{-1} \gamma^\star \Jp_{\Phat+\Pi}\,.
\end{equation}
Here $\gamma (\Phat+\Pi)^{-1} \gamma^\star$ is an $m\times m$ matrix whose $k,p$ entry when applied to $v\in\Cdotsinf(\dBCX;\Ehat)$ is defined as the limit, at $\dBCX$ from the interior of $X$,  of $D_\rho^{k-1} (\Phat+\Pi)^{-1} D_\rho^{p-1}(\delta(\rho) \times v)$. The strong transmission condition is satisfied near $\dBCX$ since $\Phat$ is a differential operator in a neighborhood of $\dBCX$ and $\Pi$ is smoothing.
So applying Corollary \ref{cor:transmission property} to $Q=D_\rho^{k-1} (\Phat+\Pi)^{-1} D_\rho^{p-1}$ we conclude that this limit exists and defines an element of $\Psi_\phi^{-m+k+p-1}(\dBCX;\Ehat)$. Since the $p,l$ entry of $\Jp_{\Phat+\Pi}$ is a differential operator of order $m+1-p-l$ if this is non-negative and equals zero otherwise it follows that $\Chat\in\Psi^*_\phi(\dBCX;\Ehat^m)$, where the order of the $k,l$ component is $k-l$.

Now fix an admissible function space $\calF$.
We apply the arguments in Section \ref{ssec:outline Calderon constr}
to $\Phat+\Pi$, using additionally that $\gamma$ and $\Jp_{\Phat+\Pi}=\Jp_\Phat$ respect $\calF$, to conclude that $\Chat$ is a projection in $\calF(\dBCX;\Ehat)^m$ to $\calB_{\Phat_X+\Pi_X;\calF}$.
Next, Proposition \ref{prop:modification B} implies
\begin{equation}
 \label{eqn:B equality}
 \calB_{\Phat_X+\Pi_X;\calF} = \calB_{\Phat_X;\calF} \,,
\end{equation}
so we conclude that $\Chat$ is an $\calF$-Calder\'on projector for $\Phat$.
The assumption of the proposition is satisfied since $\Pi$ projects to $\ker\Phat_X\cap\ker\gamma$, so the operator $\Phat_X+\Pi_X$ is a modification of $\Phat_X$, and since $\Phat$ is self-adjoint.

Finally, in view of the fact that $\Phat_X$ augments $P$, where the vector bundles $\Ehat$ and $E$ are related as in \eqref{eqn:relation E Ehat}, Proposition \ref{prop:augment B C} implies that
 \[ \calB_{P,\calF} = \pi\calB_{\Phat_X;\calF} \text{ and $C = \pi\Chat\iota$} \]
is an $\calF$-Calder\'on projector for $P$.
Now $\Chat\in\Psi^*_\phi(\dBCX;\Ehat^m)$ implies $C\in\Psi^*_\phi(\dBCX;E^m)$, with the same orders of the components. This completes the proof of Theorem \ref{thm:Calderon} (i). 
\end{proof}

\begin{proof}[Proof of Corollary \ref{cor:L2 closure of BC spaces}]
First, $\calF(\dBCX) \subset L^2_\phi(\dBCX)$ implies $\calB_{P,\calF}\subset\calH$
since $u\in\calF(X)$ implies $\gamma u\in \calF(\dBCX)^m$ (see \eqref{eqn:first.gamma}, which is proved after Definition \ref{def:admissible function space})
and $L^2_\phi(\dBCX)\supset H^k_\phi(\dBCX)$ for all $k\geq0$.

Now choose a Calder\'on projector as in Theorem \ref{thm:Calderon}.
Recall that $C=(C_{kl})_{k,l=1\dots m}$
where $C_{kl} \in \Psi_\phi^{k-l}(\dBCX;E)$. In particular, $C_{kl}$ is bounded as an operator $H^{l-1}_\phi\to H^{k-1}_\phi$, so $C$ is bounded $\calH\to\calH$. Now  $\calF^m\subset\calH$ is dense for any admissible $\calF$ since $(\Cdotsinf)^m\subset\calH$ is dense. It is a simple exercise to show that for a bounded projection on a Hilbert space the image of a dense subspace is a dense subspace of the image. Therefore, $\calB_{P,\calF}$ is dense in the range of $C$, considered as an operator on $\calH$.
\end{proof}

The proof of  Theorem \ref{thm:calderon canonical inv extension} is completely analogous to the proof in the non-singular case, see e.g.\ \cite[Lem.\ 5]{See66}. For completeness we recall the argument in our setting.

\begin{proof}[Proof of Theorem \ref{thm:calderon canonical inv extension}] Let $P$ and $\Phat$ be as in Theorem \ref{thm:calderon canonical inv extension} and $C_\pm$ be the Calder\'on projectors for $\Phat_{X^\pm}$ as constructed in the proof of Theorem \ref{thm:Calderon}, but without the augmentations and modifications.
Thus, given $U \in \Cdotsinf(\dBCX;E)^m$, let $u \in \calD^\p(\Xhat;E)$ be the solution of $\Phat u = \gamma^\star \kerningJphat U$, with $\Jphat:=\Jp_\Phat$ defined in \eqref{eqn:dfn.jump.op}. The restriction of $\Phat u$ to $\interior{\Xhat_\pm}$ is $0$, and because $\Phat^{-1}$ has the transmission property, the restriction of $u$ to $\interior{\Xhat_\pm}$ extends to a smooth function $u_\pm$ on $\Xhat_\pm$. Then $C_\pm U = \gamma_\pm u$. Moreover, $u = \Phat^{-1} \gamma^\star \Jphat U$ and since $\gamma^\star \Jphat U \in H_\phi^{1/2-m-\eps}(\Xhat;E)$ for any $\eps > 0$, this gives $u \in H_\phi^{1/2-\eps}(\Xhat;E)$. In particular, $u \in L^2_\phi(\Xhat;E)$ and so $u = u_+ + u_-$ where $u_\pm$ are extended  by zero to $\Xhat$.

Now \eqref{eqn:dfn.jump.op} and its $\Xhat^-$ counterpart read
 \begin{equation}\label{eqn:jump.op.ibp.2}
  \Phat u_\pm = \big(\Phat u_\pm\big)^0 + \gamma^\star \kerningJphat \gamma_\pm u_\pm \,.
 \end{equation}
The first term on the right vanishes since $\Phat u = \gamma^\star \Jphat U$ vanishes outside $\dBCX$ and hence so does $\Phat u_\pm$. Therefore, using $u = u_+ + u_-$ and $\gamma_\pm u_\pm = \gamma_\pm u$,
 \begin{equation}\label{eqn:jump.op.ibp.3}
  \gamma^\star \Jphat U = \Phat u = \sum_{\pm} \Phat u_\pm = \sum_{\pm} \gamma^\star \kerningJphat \gamma_\pm u_\pm = \sum_{\pm} \gamma^\star \kerningJphat \gamma_\pm u \,.
 \end{equation}
Since $\gamma^\star$ and $\Jphat$ are easily seen to be injective, this gives $U = \gamma_+ u + \gamma_- u$.
As $C_\pm U = \gamma_\pm u$ we obtain $C_+ + C_- = \Id$  on $\Cdotsinf(\dBCX;E)^m$, hence as elements of $\Psi^*_\phi(\dBCX;E)$. It also follows that $C_\Phat=C_+$.
As $C_\pm$ are bounded projections in $\calH$ with range $\calB^\pm_{P,L^2}$ we obtain the first claim and the first part of the second claim. Theorem \ref{thm:Calderon} also gives the third claim.

The normal families satisfy $N(C_+)+N(C_-)= N(\Id)=\Id$, and since
$N(C_\pm)(\mu)$ is a Calder\'on projector for $N(\Phat_{X^\pm})(\mu)$ for each $\mu$
(see Proposition \ref{prop:Calderon.normal}) the second part of the second claim also follows.
\end{proof}

\subsection{The symbol and the normal family}\label{ssec:Calderon.properties}
In this subsection we study the $\phi$-principal symbol and the normal family of the Calder\'on projector constructed above, and prove that the full $\phi$-symbol is independent of the choices made in the construction.

Given a $\phi$-elliptic operator $P \in \Diff_\phi^m(X;E,E^\prime)$ and a Calder\'on projector $C$ as in Theorem \ref{thm:Calderon}, we can follow H\"ormander's arguments (see \cite[Thm.\ 20.1.3]{Hor85}) to shed light on the $\phi$-principal symbol $\phisigma_*(C)$ of $C$.

Recall that  $\phisigma_*(C)$ is a function on (or rather a section of a homomorphism bundle over) $\phiT^*\dBCX$. To describe it we use the identification of a neighborhood of $\dBCX$ in $\Xhat$ with
$(-1,1)\times\dBCX$ from \eqref{eqn:trivialize bd}, which
induces an isomorphism
$$\phiT^*_{p} X \cong \R\times\phiT^*_{p}\dBCX;\quad \xi \mapsto (\omega,\xi')$$
for $p \in \dBCX$. Recall that $\phisigma_m(P)$ is defined on $\phiT^* X$.
For each $(p,\xi')\in \phiT^*\dBCX$ consider the ordinary differential equation
 \begin{equation}\label{eqn:symbol.ode}
  \phisigma_m(P)(p;D_t,\xi^\p) v(t) = 0 \,,\quad t\in\R\,.
 \end{equation}
Let $\Cdotinf_\pm(\R;E_{p})$ be the space of smooth functions $v: \R \rightarrow E_{p}$ for which $v(t)$ vanishes rapidly with all its derivatives as $t \to \pm\infty$.
As in \cite[Thm.\ 20.1.3]{Hor85} the $\phi$-principal symbol of $C$ is expressed in terms of solutions of \eqref{eqn:symbol.ode}:

\begin{proposition}\label{prop:Calderon.symbol} Let $(p,\xi^\p) \in \phiT^*\dBCX$ with $\xi^\p \neq 0$. The map $\phisigma_*(C)(p,\xi^\p)\in \End(E_p^m)$ is a Calder\'on projector for the differential operator $\phisigma_m(P)(p;D_t,\xi^\p)$ on $\Cdotinf_+(\R;E_{p})$.

More precisely, let $\calB_{\phisigma(P)}^\pm(p,\xi^\p)$ denote the space of boundary data at $t=0$ of solutions to \eqref{eqn:symbol.ode} with $v \in \Cdotinf_\pm(\R;E_{p})$. Then $\phisigma_*(C)(p,\xi^\p)$ is the projection in $(E_{p})^m$ with range $\calB_{\phisigma(P)}^+(p,\xi^\p)$ and kernel $\calB_{\phisigma(P)}^-(p,\xi^\p)$.
\end{proposition}

This establishes Theorem  \ref{thm:Calderon} (ii). 
The proof in \cite{Hor85}, i.e.\ in the non-singular case, is a local calculation near $\diag_{\dBCX}$, the core of which is quantifying the statement of Corollary \ref{cor:transmission property} on the level of symbols. Just as this corollary carries over to the $\phi$-case, as stated, because the local geometry near the $\phi$-diagonal $\diag_{\phi,\dBCX}$ (cf.\ Figure \ref{fig:transm.cond}) stays the same uniformly for $x\to 0$, this calculation also carries over to the singular case, so we do not repeat it here. The $\phi$-principal symbol arises since it is the symbol
in the representation, as a conormal distribution,  of the Schwartz kernel at the $\phi$-diagonal, and this symbol is defined on the dual of the normal bundle of this diagonal, which in turn is identified with $\phiT^*\dBCX$ under the projection to $\dBCX$.

\smallskip

A genuinely new feature in our singular context is the boundary symbol or normal family $N(C)$ of $C$. It turns out that the normal family of $C$ is a Calder\'on projector as well. 
Proposition \ref{prop:Calderon.normal}, below, shows part (iii) of Theorem \ref{thm:Calderon}.

\begin{proposition}\label{prop:Calderon.normal} For each $ \mu \in \R \times T^*B$, $N(C) (\mu)$ is a Calder\'on projector for $N(P)(\mu)$ on $\Cinf(F_y;E)$. More precisely, it is the Calder\'on projector resulting from the construction of Section \ref{ssec:outline Calderon constr} when starting with $N(P)(\mu)$ instead of $P$.
\end{proposition}

\begin{proof} Taking the normal family defines an algebra homomorphism and since both of $\gamma (\Phat + \Pi_{sh})^{-1} \gamma^\star$ and $\Jphat:=\Jp_\Phat$ are (matrices of) $\phi$-$\Psi$DOs on $\dBCX$ with respect to the fibration $\phi : \dBCX \to B$, we have, for each $\mu=(\tau;y,\eta)$,
 \begin{equation}\label{eqn:normal.calderon.decomp}
  N(\hat{C}) (\mu) = N(\gamma (\Phat + \Pi_{sh})^{-1} \gamma^\star) (\mu) \, N(\Jphat) (\mu) \,.
 \end{equation}
Note that all three operators in \eqref{eqn:normal.calderon.decomp} operate on sections over $\partial F_y$, which is locally given by $z_1=0$.
Recall from Appendix \ref{sec:app phi ops} the definition of the normal family for $\phi$-pseudodifferential operators, see Equation \eqref{eqn:normal op osc testing}.
We use this with $X$ replaced by $\dBCX$ and $F_y$ replaced by $\partial F_y$
to examine $N(\gamma(\Phat + \Pi_{sh})^{-1}\gamma^\star)(\mu)$ more closely.
Note that the function $g$ in  \eqref{eqn:normal op osc testing}  is constant along the fibre $\partial F_y$ and that $\gamma$, $\gamma^\star$ act in the fibre direction $z_1$ only. Therefore, $\gamma$ and $\gamma^\star$ commute with multiplication by $e^{\pm i g}$.
Now define $\gamma_y$ and $\gamma_y^\star$ analogously to $\gamma$ and $\gamma^\star$, but fibre-wise for $\partial F_y \subset F_y$. Then if $V \in C^\infty(\dsBCX;\Ehat)^m$ and $\tilde V$ denotes a smooth extension of $V$ to $\dBCX$, we see that $\widetilde{\gamma_y^\star V}:=\gamma^\star\tilde V$
restricts to $\gamma_y^\star V$ at $\dsBCX$. Therefore,
 \begin{align}
  N(\gamma (\Phat &+ \Pi_{sh})^{-1} \gamma^\star) (\mu) V
    = \big[ e^{-ig} \gamma (\Phat + \Pi_{sh})^{-1} \gamma^\star (e^{ig}\tilde V)\big]_{|\partial F_y} \\
   &= \big[\gamma e^{-ig} (\Phat + \Pi_{sh})^{-1} (e^{ig}\gamma^\star \tilde V)\big]_{|\partial F_y}
    = \gamma_y \big[ e^{-ig} (\Phat + \Pi_{sh})^{-1} (e^{ig} \widetilde{\gamma_y^\star V})\big]_{|\partial F_y} \\
   &= \gamma_y N((\Phat + \Pi_{sh})^{-1})(\mu)\gamma_y^\star \, V \,.
  \end{align}

As, locally near $\dsBCX$, the construction of the jump terms $\Jphat = \Jp_{\Phat} = \Jp_{\Phat + \Pi_{sh}}$ involved the fibre direction $z_1$ only, by a similar argument we see that $N(\Jphat) (\mu)$ are the jump terms for $N(\Phat + \Pi_{sh})$, that is $\Jp_{N(\Phat+\Pi_{sh})}$, compare \eqref{eqn:dfn.jump.op}. Since $N((\Phat + \Pi_{sh})^{-1}) = (N(\Phat+\Pi_{sh}))^{-1}$ and $N(\Pi_{sh})=0$ (as the Schwartz kernel of $\Pi_{sh}$ vanishes rapidly at $(\dsX)^2$), this shows that
 \[ N(\hat{C}) = \gamma_y N((\Phat+\Pi_{sh})^{-1}) \gamma_y^\star \Jp_{N(\Phat+\Pi_{sh})}
               = \gamma_y (N(\Phat))^{-1} \gamma_y^\star \Jp_{N(\Phat)} \]
is, by the same arguments that lead to Theorem \ref{thm:Calderon}, a Calder\'on projector for $N(\Phat)$ or more precisely for $N(\Phat)_{|X} = N(\Phat_{|X})$.

Now, using \eqref{eqn:normal op osc testing} again, we see that taking the normal family commutes with the bundle maps $\pi$ and $\iota$, i.e., $\pi N(\Chat) \iota = N(\pi \Chat \iota) = N(C)$ by Proposition \ref{prop:augment B C} and $\pi N(\Phat_{|X}) \iota = N(\pi \Phat_{|X} \iota) = N(P)$ by Theorem \ref{thm:construction inv ext}.
Also, taking inverses and products and taking $\Jp$ is compatible with the $\pi$-$\iota$-restriction, as is most easily seen using the direct summand characterization of augmentation, Equation \eqref{eqn:augment op direct sum}.
By Proposition \ref{prop:augment B C} again, we obtain the claim and the proof is complete.
\end{proof}

\begin{proof} [Proof of Proposition \ref{prop:Calderon.unique}]
As in the proof of Proposition \ref{prop:Calderon.normal} we use that $C=\pi \Chat\iota$ where $\Chat = \gamma (\Phat + \Pi_{sh})^{-1} \gamma^\star\Jphat$. For $p\in\Xhat$ denote by $\sigmafull(Q)(p)$ the full $\phi$-symbol of a $\phi$-$\Psi$DO $Q$ at $p$. Using a choice of coordinates near $p$ this is defined as a function on ${}^\phi T^*_p\Xhat$ modulo Schwartz functions. By the standard parametrix construction, for any $k$ the $k$-jet of $\sigmafull((\Phat+\Pi_{sh})^{-1})$ at $p$ is determined by the infinity-jet of $\sigmafull(\Phat+\Pi_{sh})=\sigmafull(\Phat)$ at $p$.
Since the same holds for  $\sigmafull(\Jphat)$, it follows for $p\in\dBCX$ that $\sigmafull(\Chat)(p)$ is determined by the infinity-jet of $\sigmafull(\Phat)$, hence of $\sigmafull(\Phat_X)$, at $p$.
Applying $C=\pi \Chat\iota$ and $\pi\Phat_X\iota=P$ we obtain as at the end of the proof of Proposition \ref{prop:Calderon.normal} that $\sigmafull(C)(p)$
is determined by the infinity-jet of $\sigmafull(P)$ at $p$.
This completes the proof.
Note that instead of the inverse of $\Phat+\Pi_{sh}$ a parametrix modulo $\Psi_\phi^{-\infty}$
in a neighborhood of $\dBCX$ suffices to fix $C$ modulo $\Psi_\phi^{-\infty}(\dX)$.
\end{proof}

\subsection{The orthogonal Calder\'on projector}\label{ssec:orth proj}
We now prove Theorem \ref{thm:orthogonal}. It follows from the following proposition.
\begin{proposition}\label{thm:orth.Calderon} Let $C$ be a Calder\'on projector for $P \in \Diff_\phi^1(X;E,E^\p)$ as constructed in Theorem \ref{thm:Calderon}. Denote the $L^2_\phi$-adjoint of $C$ by $C^\star$.

Then $\Id + C - C^\star$ is a fully elliptic and invertible $\phi$-$\Psi$DO and
 \begin{equation}\label{eqn:orth.calderon}
  C_o = C ( \Id + C - C^\star)^{-1} : L^2_\phi(\dBCX;E) \longrightarrow \calB_{P,L^2}
 \end{equation}
is the orthogonal projection (with respect to the $L^2_\phi$-scalar product).

We have $C_o \in \Psi_\phi^0(\dBCX;E)$, with the $\phi$-principal symbol and the normal family of $C_o$ being the orthogonalizations (obtained as in \eqref{eqn:orth.calderon}) of the $\phi$-principal symbol respectively the normal family of $C$.
\end{proposition}

The main tool in proving Proposition \ref{thm:orth.Calderon} is part b) of Lemma 3.5 of \cite{BLZ09} dealing with orthogonalizations of bounded projections: it states that whenever $T : H \to H$ is a bounded projection in a Hilbert space $H$ and $T^\star$ is its Hilbert space adjoint, then $\Id + T - T^\star$ is invertible and
 \begin{equation}\label{eqn:orth.formula}
  T_o = T (\Id + T - T^\star)^{-1} : H \longrightarrow H
 \end{equation}
is the orthogonal projection with range $\rg T$.

\begin{proof}[Proof of Proposition \ref{thm:orth.Calderon}]
By construction $C$ is a projection, and since $C \in \Psi_\phi^0(\dBCX;E)$, it is bounded on $L^2_\phi(\dBCX;E)$. The operators $\phisigma_0(C)(p^\prime,\xi^\prime)$ and $N(C)(\mu)$ are bounded projections as well, by Propositions \ref{prop:Calderon.symbol} and \ref{prop:Calderon.normal}.

We now apply part b) of Lemma 3.5 of \cite{BLZ09} separately to $C$, its $\phi$-principal symbol and its normal family. In the following, we omit mentioning the variables $(p^\prime,\xi^\prime)$ and $(\tau;y,\eta)$ for the $\phi$-principal symbol and normal family and denote by a star $\star$ the appropriate Hilbert space adjoint. For $C$ and $N(C)$, this will be the $L^2_\phi$- respectively $L^2$-adjoint and for $\phisigma_0(C)$, this will be the adjoint of an endomorphism of the hermitian bundle $E$.

As the $\phi$-principal symbol- and normal family-maps are $\star$-algebra homomorphisms (see \cite{MM98}), we have
\begin{align}
   \phisigma_0(\Id+C-C^\star)
     &= \Id + \phisigma_0(C) - \phisigma_0(C)^\star \,, \\
   N(\Id + C - C^\star)
     &= \Id + N(C) - N(C)^\star \,,
 \end{align}
so by Lemma 3.5 (loc.cit.)\ these operators, as well as  $\Id + C - C^\star$,
are invertible. Thus, $\Id + C - C^\star$ is fully elliptic and invertible on $L^2_\phi(\dBCX;E)$, whence its inverse is in $\Psi_\phi^0(\dBCX;E)$ again. But then
 \begin{equation}\label{eqn:orth.Calderon.red}
  C_o = C \big( \Id + C - C^\star\big)^{-1} \in \Psi_\phi^0(\dBCX;E)
 \end{equation}
with $\phi$-principal symbol
 \begin{equation}\label{eqn:orth.Calderon.red.symb}
  \phisigma_0(C_o)
     = \phisigma_0(C) \big( \Id + \phisigma_0(C) - \phisigma_0(C)^\star\big)^{-1}
 \end{equation}
and normal family
 \begin{equation}\label{eqn:orth.Calderon.red.normal}
  N(C_o)
    = N(C) \big( \Id + N(C) - N(C)^\star\big)^{-1} \,,
 \end{equation}
and by the lemma again, \eqref{eqn:orth.Calderon.red}, \eqref{eqn:orth.Calderon.red.symb} and \eqref{eqn:orth.Calderon.red.normal} are the orthogonalizations of $C$, $\phisigma_0(C)$ and $N(C)$, respectively. Clearly, $C_o$ has the same range as $C$, i.e.\ $\calB_{P,L^2}$.
\end{proof}

\begin{remark}\label{rem:orthogonal}
\mbox{}
\begin{enumerate}
 \item
 We chose to formulate the theorem and proposition on orthogonal projections only for $m=1$ since for $m>1$ one needs to introduce scalar products in the Sobolev spaces $\calH$ resp.\ $H^k_\phi(\dBCX;E)$, which seems (even) less natural. However, it is easy to generalize the statements and proofs to this case, using order reduction operators in the $\phi$-calculus.
\item
 In general, the $\phi$-principal symbols of $C$ and $C_o$ will be different.
 Heuristically this is clear since the range and kernel spaces $\calB_{\phisigma(P)}^\pm(p^\p,\xi^\p) \subset (E_{p'})^m$ of $\phisigma(C)(p',\xi')$ (see Proposition \ref{prop:Calderon.symbol})
 cannot be expected to be orthogonal to each other for each value of $\xi'$, as they would have to be if $\phisigma(C)=\phisigma(C_o)$.
For example, if $P$ is the Laplacian of a $\phi$-metric then a simple calculation shows that this is indeed not the case.
\end{enumerate}
\end{remark}

\appendix
\section{Basics on Manifolds with Corners and Blow-Ups} \label{sec:app mwc}
We give a quick summary of basic notions on manifolds with corners. Details can be found in \cite{Mel:APSIT}, \cite{Mel:DAMWC} or in the introductory text \cite{Gri:BBC}.

A \textbf{manifold with corners} of dimension $n$, denoted $X$ in the sequel, is defined like a manifold except that local charts are defined on open subsets of model spaces $\R^n_k:=\R^{n-k}\times\R_+^k$ for various $k\in\{0,\dots,n\}$ where $\R_+=[0,\infty)$, and an additional global condition is satisfied, see below. To define smoothness of transition maps or of maps between manifolds with corners we say that a map
 \[ \text{$U_1\to U_2$, where $U_i\subset\R^{n_i}_{k_i}$,} \]
is \textbf{smooth} if it extends to a smooth map
 \[ \text{$\Utilde_1\to \Utilde_2$, where $\Utilde_i\subset\R^{n_i}$ are open and $U_i=\tilde U_i\cap \R^{n_i}_{k_i}$.} \]
If $p\in X$ then there is a unique $k$, called the \textbf{codimension} of $p$, so that there is a coordinate system (inverse of a chart) mapping $p$ to $0\in\R^n_k$. The coordinates are then sometimes called \textbf{adapted} to $X$, and are often denoted $x=(x_1,\dots,x_k)$, $y=(y_1,\dots,y_{n-k})$ where $x_i\geq0$ and $y_j\in\R$ for all $i,j$.

A \textbf{face} of $X$ of codimension $k$ is the closure of a connected component of the set of points of codimension $k$. A \textbf{boundary hypersurface} is a face of codimension one. The global condition on a manifold with corners is that boundary hypersurfaces be embedded (rather than immersed) submanifolds. Equivalently, for each boundary hypersurface $H$ there is \textbf{boundary defining function} $\rho$, i.e.\ a smooth function $\rho:X\to\R_+$ satisfying $\rho^{-1}(0)=H$ and $d\rho_{|p}\neq0$ for all $p\in H$.

A \textbf{p-submanifold} (where p is for product) of  $X$ is a subset $Y$ so that for each $p\in Y$ there is an adapted coordinate system on $X$ in which $Y$ is locally a coordinate subspace. That is, adapted coordinates $z=(x,y)$ can be chosen and  regrouped as $(z',z'')$ so that $Y=\{z''=0\}$.
Also, $Y$ is called a \textbf{boundary p-submanifold} if $Y\subset\dX$, otherwise it is an \textbf{interior p-submanifold}. In the latter case only $y$ variables occur among the $z''$ variables.
For example, faces of $X$ are boundary p-submanifolds.

If $Y\subset X$ is a closed p-submanifold then the \textbf{blow-up} of $X$ in $Y$ is a new manifold with corners, denoted $[X,Y]$, together with a smooth map $\beta:[X,Y]\to X$, called \textbf{blow-down map}, which restricts to a diffeomorphism $[X,Y]\setminus\ff \to X\setminus Y$, where $\ff:=\beta^{-1}(Y)$ is called the \textbf{front face}, and so that
near any $p\in Y$ with coordinates $(z',z'')$ as above the map $\beta$ is locally near $\beta^{-1}(p)$ modelled by the polar coordinates map in the $z''$-coordinates, i.e.\ if $z'\in\R^{n'}_{k'}$, $z''\in\R^{n''}_{k''}$ then locally $\beta: \R^{n'}_{k'}\times \R_+\times \Sph^{n''-1}_{k''}\to\R^{n'}_{k'}\times\R^{n''}_{k''}, (z',r,\omega)\mapsto (z',r\omega)$ where $\Sph^{n''-1}_{k''}\subset\R^{n''}_{k''}$ is the unit sphere. Locally the front face is $\R^{n'}_{k'}\times\Sph^{n''-1}_{k''}$ and has local boundary defining function $r$. In practice, it is better to use \textbf{projective coordinates}. Examples are given in Appendix \ref{sec:app phi ops} for the double spaces $X^2_b$ and $X^2_\phi$.
Blow-ups can also be iterated, i.e.\ if $Z$ is a p-submanifold of $[X,Y]$ then one can form $[[X,Y],Z]$ etc.

If $Z$ is a connected subset of $X$ then the \textbf{lift} of $Z$ under the blow-up of a p-submanifold $Y\subset X$, denoted $\beta^*(Z)$, is defined as $\beta^{-1}(Z)$ if $Z\subset Y$ and as the closure of $\beta^{-1}(Z\setminus Y)$ otherwise. If $Z$ is a p-submanifold meeting $Y$ \textbf{cleanly} (i.e.\ so that for each $p\in Y\cap Z$ there is an adapted coordinate system in which both $Y$ and $Z$ are coordinate subspaces) then $\beta^*(Z)$ is a p-submanifold of $[X,Y]$. However, also subsets $Z$ which are not p-submanifolds can become such after (possibly iterated) blow-up, and then we say that $Z$ is \textbf{resolved} by the (iterated) blow-up. An important example is the diagonal $\{x=x'\}$ in $\R^2_+$, which is not a p-submanifold but is resolved by blowing up the origin.
A p-submanifold $Y$ meets any face of $X$ cleanly, and the boundary hypersurfaces of $[X,Y]$ are $\ff$ and the lifts of the boundary hypersurfaces of $X$.

From now on we assume for simplicity that $X$ is compact.
The space of all smooth vector fields on $X$ which at each $p\in X$ are tangent to all boundary hypersurfaces containing $p$ is denoted $\calV_b(X)$.
Interpreting vector fields as first order differential operators and taking finite sums of smooth functions and compositions $V_1\circ\dots\circ V_l$ with all $V_i\in\calV_b(X)$, for $l\leq m$, we obtain the space of \textbf{b-differential operators} of order at most $m$, denoted $\Diff_b^m(X)$. Also $\Diff_b^*(X) := \bigcup_m \Diff_b^m(X)$. In adapted coordinates these are combinations of expressions $x_i\partial_{x_i}$, $\partial_{y_j}$ with smooth coefficients.

We now define various function spaces on a manifold with boundary $X$.\footnote{The definitions can be extended to  manifolds with corners, but we don't need this.} First, $\Cinf(X)$ denotes  the space of functions $X\to\C$ which are smooth up to the boundary.
Let $x$ be a boundary defining function.
The following spaces consist of functions only defined and smooth on the interior of $X$ but having a certain prescribed behavior near the boundary:
\begin{itemize}
 \item The space of \textbf{functions conormal to the boundary} of order $a\in\R$:
$$\calA^a(X)= \{u\in\Cinf(\interior{X}): \Diff_b^*(X)u \subset x^a L^\infty(X)\}\,.$$
There is also the $L^2$-variant of $\calA^a(X)$, where $L^\infty$ is replaced by $L^2$ in the definition. Also, the space of all conormal functions is denoted
$$\calA(X) = \bigcup_{a\in\R}\calA^a(X)\,.$$
\item Spaces of functions \textbf{polyhomogeneous} at the boundary:
$$\calA^\calE(X) = \{u\in\calA(X):\, u\sim \sum_{(z,k)\in \calE} a_{z,k} \,x^z\log^k\!x\}
$$
where $\calE$ is an \textbf{index set}, i.e.\ a discrete subset of $\C\times\N_0$ satisfying certain additional conditions,  and $a_{z,k}\in\Cinf(X)$.
The symbol `$\sim$' means that the difference of $u$ and the finite sum over $\Re z<N$ is in $\calA^N(X)$, for each $N$.
If $\calE=\N_0\times\{0\}$ then $\calA^\calE(X)=\Cinf(X)$. If the boundary is disconnected then each component can have its own index set (resp. order $a$ for $\calA^a(X)$).
\end{itemize}
If $X$ is a manifold with \textBC-boundary (see Definition \ref{def:bv-mfd}) whose singular boundary $\dsX$ is a disjoint union of boundary hypersurfaces
then we denote by $$\calA^a_\s(X),\ \calA_\s(X),\ \calA^\calE_\s(X)$$ the corresponding spaces where we assume the respective behavior only to occur at $\dsX$, smoothly up to $\dBCX$. Alternatively, elements of $\calA^a_\s(X)$ are restrictions to $X$ of elements of $\calA^a(\Xhat)$ where $\Xhat$ is the \textBC-double of $X$, cf.\ Section \ref{ssec:extensions}, and similarly in the other cases.

We denote by $L^2_\phi(X)$ the $L^2$-space with respect to the volume form of a $\phi$-metric, and define $\phi$-Sobolev spaces $H^k_\phi(X):=\{u:\,\Diff^k_\phi(X)u\subset L^2_\phi(X)\}$ for $k\in\N_0$.

For the definition of pseudodifferential operators we  need conormal distributions. A distribution $u$ on a manifold with corners $Z$ is \textbf{classical conormal} of order $m\in\R$ with respect to an interior p-submanifold $Y$ if it is smooth on $Z\setminus Y$ and near any point of $Y$ and in any adapted coordinate system $x,y=(y',y'')$ for $Z$ in terms of which $Y=\{y''=0\}$ locally,
\begin{equation}
 \label{eqn:def conormal}
u(x,y',y'') = \int e^{iy''\eta''}a(x,y';\eta'')\,d\eta''
\end{equation}
for a classical symbol $a$ of order $\mu=m+\frac14\dim Z - \frac12\codim Y$.
The space of these distributions is denoted $I^m_\cl(Z,Y)$.
We only need the case $\dim Z=2\dim Y$, then $\mu=m$.
Here classical means that $a$ has a symbol expansion
$a\sim \sum_{j=0}^\infty a_j$ where for each $j$
\begin{equation}
\label{eqn:pos homogeneous}
a_j (x,y';\lambda\eta'') = \lambda^{\mu-j} a_j(x,y';\eta'')
\end{equation}
for all $\lambda>0$ and all $x,y'$ and $\eta''\neq0$.
The $a_j$ are uniquely determined by $u$.

\section{Background on Fibred Cusp Operators} \label{sec:app phi ops}

In Section \ref{ssec:phi-geom} we defined $\phi$-manifolds, $\phi$-metrics, $\phi$-vector fields, $\phi$-differential operators and their $\phi$-symbol and normal family.
Here we describe the $\phi$-pseudodifferential calculus introduced in \cite{MM98}. We assume that $M$ is a $\phi$-manifold without \textBC-bondary, i.e.\ simply a manifold with boundary (here we write $\dM$ rather than $\dsM$), which at first we assume to be connected, equipped with a fibration $\phi:\dM\to B$ and a boundary defining function $x$.

We want to define a pseudodifferential calculus, i.e.\ a set of operators closed under composition, which contains $\Diff_\phi^*(M)$, as defined in Section \ref{ssec:phi-geom}, as well as parametrices of fully elliptic elements of this space. By the general philosophy on singular pseudodifferential calculi introduced by R. Melrose,
such a calculus is defined by a set of Schwartz kernels, which are distributions on the (interior of the) double space $M^2=M\times M$, and whose boundary behavior is restricted by requiring that their pull-backs to a suitable blow-up, $M^2_\phi$, of $M^2$ satisfy certain smoothness and vanishing conditions at the boundary hypersurfaces of $M^2_\phi$ and have a conormal singularity at the diagonal, uniformly up to the boundary.

The $\phi$-double space, $M^2_\phi$, is defined as follows: First, blow up $(\dM)^2\subset M^2$. This yields the b-double space with blow-down map
$$\beta_b: M^2_b := [M^2,(\dM)^2] \to M^2 \,. $$
Its front face is denoted by $\bbf$, and the interior of $\bbf$ is naturally diffeomorphic to $(0,\infty)\times(\dM)^2$ where the first coordinate is $t=\frac{x'}x$, with $x,x'$ the pull-backs of the boundary defining function $x$ on $M$ to the first and second factor in $M^2$.
Let $\diag_M\subset M^2$ be the diagonal and $\diag_{b}$ be its lift to $M^2_b$. It meets the boundary of $M^2_b$ in the interior of $\bbf$, in the set $\{1\}\times\diag_\dM$. The larger interior submanifold $D_\phi:=\{1\}\times\{(p,p')\in(\dM)^2:\phi(p)=\phi(p')\}$ of $\bbf$ is called the fibre diagonal.
We blow this up and define the $\phi$-double space
$$ \beta_{D_\phi}: M^2_\phi := [M^2_b,D_\phi]\to M^2_b,\quad  \beta_\phi:= \beta_b\circ\beta_{D_\phi}:M^2_\phi\to M^2\,.$$
The front face created by this blow-up is denoted $\phif$.
The diagonal $\diag_{b}$ lifts to a p-submanifold $\diag_\phi$ of $M^2_\phi$ which meets the boundary in the interior of $\phif$.

We describe these spaces locally, using adapted local coordinates $x,y,z$ (see Section \ref{ssec:phi-geom}; recall that $\phi(x,y,z)=(x,y)$). Starting from coordinates $(x,y,z; x',y',z')$ on $M^2$
 we have $(\dM)^2= \{x=x'=0\}$, so coordinates near interior points of $\bbf$ are
 $$ t=\frac {x'}{x}, x, y,z,y',z'\,,$$
 with $x$ defining $\bbf$ there.
The diagonal in $M^2$ is $\{x=x',y=y',z=z'\}$ and lifts to $\{t=1,y=y',z=z'\}$
 The fibre diagonal $D_\phi$ is locally $\{t=1, x=0, y=y'\}$. So projective coordinates near the interior of $\phif$ are
\begin{equation}
\label{eqn:phi double coords}
 T := \frac{t-1}{x}, x, Y:=\frac{y-y'}{x}, y, z,z'
\end{equation}
 with $x$ defining $\phif$ there.
 In these coordinates the diagonal $\diag_\phi$ is given by $\{T=0, Y=0, z=z'\}$, and $|(T,Y)|\to\infty$ corresponds to the boundary of $\phif$, which is its intersection with the lift of $\bbf$.
The fundamental reason for considering the space $M^2_\phi$ is that $\phi$-vector fields on $M$, when considered as vector fields on $M^2$ in the $x,y,z$ variables, i.e.\ pulled back from the left factor, lift under $\beta_\phi$ to smooth vector fields on $M^2_\phi$ that span a rank $\dim M$-bundle and at any $\gamma\in \diag_\phi$ span a subspace of $T_\gamma M^2_\phi$ transversal to $T_\gamma\diag_\phi$.

We consider operators $P$ acting on functions on $\interior{M}=M\setminus\dM$ which are given by Schwartz kernels $K_P$, which are distributions on $(\interior M)^2=\interior{M^2}$, in the sense that
\begin{equation}
 \label{eqn:schwartz kernel}
 (Pu)(p) = \int_{\interior M} K_P(p,p') u(p')\, \nu(p')
\end{equation}
where $\nu$ is some fixed density on $\interior M$. We choose (and fix once and for all) for $\nu$ a smooth positive $\phi$-density, i.e.\ locally $\nu = a \frac{dx}{x^2}\frac{dy}{x^b}dz$ with $a>0$ smooth up to the boundary $x=0$ and $b = \dim B$.
For example, $\nu$ could be the volume density of a $\phi$-metric. The reason for this choice is that for $P=\Id$ (and, say, $a=1$), we have $K_P=(x')^{b+2}\delta(x-x')\delta(y-y')\delta(z-z')$, which in coordinates \eqref{eqn:phi double coords} is $\delta(T)\delta(Y)\delta(z-z')$, and since this has no $x$-factor, it extends from the interior of $M^2_\phi$  to a distribution on $M^2_\phi$ as a smooth non-vanishing delta-function for the submanifold $\diag_\phi$.
The lifting property of $\phi$-vector fields mentioned above then implies that the Schwartz kernels for $P\in\Diff^m_\phi(M)$ lift to $\interior{M^2_\phi}$ and extend to $M^2_\phi$ to be delta-functions of order at most $m$ for $\diag_\phi$.
Explicitly, if $P$ is given in coordinates by \eqref{eqn:phi operator} and $\nu=\frac{dx}{x^2}\frac{dy}{x^b}dz$ then
\begin{equation}
 \label{eqn:kernel of phidiff op}
 K_{P} = \sum_{k+|\alpha|+|\beta|\leq m} a_{k,\alpha,\beta}(0,y,z) D_T^k\delta(T) D_Y^\alpha \delta(Y) D_z^\beta\delta(z-z') + O(x)
\end{equation}
We define $\phi$-pseudodifferential operators by replacing delta-functions by the larger space of classical conormal distributions, defined in Appendix \ref{sec:app mwc}:

\begin{definition}
Let $M$ be a $\phi$-manifold and $m\in\R$. The space $\Psi_\phi^m(M)$ is defined as the set of operators \eqref{eqn:schwartz kernel} whose Schwartz kernels $K_P$ lift to $M^2_\phi$ to elements of
$$ I^m_\cl(M^2_\phi, \diag_\phi) $$
that vanish to infinite order at all boundary hypersurfaces of $M^2_\phi$ except $\phif$.
\end{definition}

There is also a more general definition without $\cl$, but all operators occurring in this paper are classical. By the remarks before the definition we have $\Diff_\phi^*(M) \subset \Psi_\phi^*(M)$. The definition of the $\phi$-principal symbol extends in a straight-forward way using the local representation \eqref{eqn:def conormal}. The definition of the normal family of $P\in\Diff_\phi^m(M)$ given in \eqref{eqn:def normal op diff} does not extend directly to $\Psi_\phi^*(M)$. However, reinterpreting this formula in terms of \lq oscillatory testing\rq\ allows to extend it as follows: for $\tau\in\R$, $\eta\in\R^b$ and $y_0\in B$ let $g(x,y)=-\frac\tau x + \frac\eta x (y-y_0)$ in coordinates on $B$ near $y_0$.
Then $x^2D_x e^{ig} = \left(\tau-\eta(y-y_0)\right)e^{ig}$ and $xD_y e^{ig} = \eta e^{ig}$, and this implies that for $P$ as in \eqref{eqn:phi operator} and $u\in\Cinf(M)$ we have
\begin{equation}
 \label{eqn:normal op osc testing}
 \left[e^{-ig}P(e^{ig}u)\right]_{|F_{y_0}} = N(P)(\tau;y_0,\eta)(u_{|F_{y_0}})
\end{equation}
as functions on the fibre $F_{y_0}$, since $F_{y_0}$ is given by $x=0, y=y_0$.
It can be shown that the left hand side is well-defined and smooth for $P\in \Psi^*_\phi(M)$ and only depends on $u_{|F_{y_0}}$, and that $N(P)$ so defined is a (standard) pseudodifferential operator with parameter $(\tau,
\eta)$ on $F_{y_0}$, and varies smoothly in $y_0$. Also, the definition makes sense invariantly when considering $\eta\in T_{y_0}B$. By \eqref{eqn:def normal op diff} the normal family for $P\in\Diff^m_\phi(M)$ vanishes if and only if $P\in x\Diff^m_\phi(M)$, and an analogous statement holds for $P\in\Psi^* _\phi(M)$.
A short calculation shows that the Schwartz kernel $K_{N(P)}(\tau;y,\eta;z,z')$ of $N(P)(\tau;y,\eta)$ is the $(T,Y)\to(\tau,\eta)$ Fourier transform of the restriction of $K_P$ to $\phif$, when writing $K_P$ in coordinates \eqref{eqn:phi double coords}.

As for $\phi$-differential operators, a $\phi$-pseudodifferential operator is called $\phi$-elliptic if its $\phi$-principal symbol is invertible outside the zero section and fully elliptic if in addition its normal family is invertible for all $\tau,y,\eta$.

The main facts about the $\phi$-pseudodifferential calculus are:
\begin{enumerate}
\item $\Psi^*_\phi(M):= \bigoplus_{m\in\R}\Psi^m_\phi(M)$ is an $\R$-graded $\star$-algebra, i.e.\ a vector space and closed under adjoints and composition, with orders adding under composition.
\item The $\phi$-principal symbol $P\mapsto \phisigma(P)$ and the normal family $P\mapsto N(P)$ are $\star$-algebra homomorphisms, i.e.\ they are linear and respect composition and the involution $\star$.
\item Operators in $\Psi^m_\phi(M)$ are bounded $H^s_\phi(M)\to H^{s-m}_\phi(M)$ for all $s$ and map each of the spaces $\calA^a(M)$, $\calA^{\calE}(M)$ into itself, for any $a\in\R$ and index set $\calE$.
\item An element of $x^a\Psi^m_\phi(M)$ is a compact operator in $L^2_\phi(M)$ if and only if $m<0$ and $a>0$.
\item An operator $P\in \Psi^m_\phi(M)$ is $\phi$-elliptic if and only if it has a parametrix with remainders in $\Psi^{-\infty}_\phi(M)$.
\item An operator $P\in \Psi^m_\phi(M)$ is fully elliptic if and only if it has a parametrix with remainders in $x^\infty\Psi^{-\infty}_\phi(M)$, if and only if it is Fredholm as a map $H^s_\phi(M)\to H^{s-m}_\phi(M)$ for any $s$.
In particular, $\ker P\subset \Cdotsinf(M)$ in this case.
\item If $P\in\Psi^m_\phi(M)$ is invertible as an operator $H^{s}_\phi(M)\to H^{s-m}_\phi(M)$ for some $s$ then its inverse is in $\Psi^{-m}_\phi(M)$.

\end{enumerate}

We make some additional remarks. Above we assumed that the boundary $\dM$ is connected. More generally, if $\dM$ has components $H_1,\dots,H_r$ (each one with its own base and fibre) then $\Psi^m_\phi(M)$ is defined in the same way, except that in the definition only submanifolds meeting the diagonal are blown up:
$M^2_b$ is defined by blowing up  $H_i^2$ for $i=1,\dots,r$ (but not $H_i\times H_j$ with $i\neq j$), and similarly for $M^2_\phi$. This yields a disjoint collection of $\phi$-faces $\phif_i$, and a normal family for each boundary component $H_i$.

The $\phi$-double space $M^2_\phi$ is defined also when $M$ is a $\phi$-manifold with \textBC-boundary. Note that the lifted diagonal $\diag_b$ is not a p-submanifold in $M^2_b$. However, the fibre diagonal $D_\phi$ is a p-submanifold (essentially since the \textBC-boundary arises from fibres having boundary, not the base), so the second blow-up used to define $M^2_\phi$ is still defined. Of course $\phi$-pseudodifferential operators on $M$ are not defined unless the \textBC-boundary is  empty.


\bibliographystyle{amsalpha}
\bibliography{D-N-cusp}

\end{document}